\documentclass[preprint,12pt]{elsarticle}

\usepackage{amsmath,amsfonts}
\usepackage{amssymb, amsxtra, mathrsfs}
\usepackage{algorithmic}
\usepackage{algorithm}
\usepackage{array}
\usepackage[caption=false,font=normalsize,labelfont=sf,textfont=sf]{subfig}
\usepackage{textcomp}
\usepackage{stfloats}
\usepackage{float}
\usepackage{url}
\usepackage{verbatim}
\usepackage{graphicx,color}
\usepackage{natbib}
\usepackage{multirow}

\usepackage[left=1in, right=1in, top=1in, bottom=1in]{geometry}
\usepackage{setspace}
\journal{arXiv}

\begin{document}

\begin{frontmatter}



\title{Locational Marginal Pricing of Energy in Pipeline Transport of Natural Gas and Hydrogen with Carbon Offset Incentives}


\author[inst1]{Mo Sodwatana}

\affiliation[inst1]{organization={Department of Energy Science and Engineering, Stanford University},
            city={Stanford},
            state={CA},
            country={USA}}

\author[inst2]{Saif R. Kazi}
\author[inst2]{Kaarthik Sundar}
\author[inst1]{Adam Brandt}
\author[inst2]{Anatoly Zlotnik}

\affiliation[inst2]{organization={Los Alamos National Laboratory},
            city={Los Alamos},
            state={NM},
            country={USA}}

\begin{abstract}
We propose an optimization formulation for locational pricing of energy transported through a pipeline network that carries mixtures of natural gas and hydrogen from distributed sources to consumers.  The objective includes the economic value provided by the pipeline to consumers of energy and suppliers of natural gas and green hydrogen, as well as incentives to lower carbon emissions by consuming the latter instead of the former.  The optimization is subject to the physics of gas flow and mixing in the pipeline network as well as engineering limits.  In addition to formulating this mathematical program, we synthesize the Lagrangian and derive analytical expressions for the dual variables.  We propose that the dual solution can be used to derive locational marginal prices of natural gas, hydrogen, and energy, as well as the decarbonization premium paid by consumers that receive hydrogen.  We derive several properties of solutions obtained using the proposed market mechanism, and demonstrate them using case studies for standard 8-node and 40-node pipeline test networks.  Finally, we show that optimization-based analysis of the type proposed here is critical for making sound decisions about economic policy and infrastructure expansion for blending green hydrogen into existing natural gas pipelines.
\end{abstract}


\end{frontmatter}


\section{Introduction}

Climate change, caused by uncontrolled greenhouse gas (GHG) emissions during the past century \cite{crowley2000causes}, is one of the most pressing global challenges today \cite{mendelsohn2004impact}. The use of fossil fuels for energy production and industrial processes is the major contributor to GHG emissions, and replacing fossil fuels with cleaner alternatives is critical to meeting the emissions reduction targets set forth \cite{owusu2016review}. The transition away from an energy system that is fossil fuel-intensive to one that is low-carbon will require a mix of renewable energy generation distributed throughout the grid, large-scale energy storage technologies, and carbon-free chemical energy carriers.

Green hydrogen refers to hydrogen gas produced via electrolysis using electricity from renewable sources. Green hydrogen is considered a promising alternative to fossil fuels, because it can serve as an alternative energy carrier and feedstock in hard-to-abate industries such as the petrochemical sectors, the cement and steel-making industries, and in heavy-duty transport \cite{yang2022}. Additionally, green hydrogen can act as a carrier or storage medium for otherwise unusable renewable generation, because so-called curtailed generation is not usable due to location and/or timing of generation. In this way, hydrogen facilities can also be integrated into power systems and offer an alternative form of long-duration or seasonal storage \cite{egeland2021integrateh2,andersson2019storageh2}.  In addition to its use in fuel cells, hydrogen can be burned to produce heat or drive turbines, so it can be injected into existing gas pipelines so that the blend of hydrogen and natural gas can be consumed in end-use appliances and furnaces \cite{topolski2022blending}. Challenges with blending of hydrogen include leakage, material degradation, and embrittlement of steel, as well as changes to distribution system pressures and possible violation of gas quality standards (e.g., lowered volumetric heating value) \cite{melaina2013blending,raju2022}. In the case that engineering issues are addressed, conceptual challenges remain with respect to characterizing the economic impacts of pipeline hydrogen blending  \cite{zlotnik2023review}.

Given the differences in the physical and chemical characteristics of hydrogen and natural gas, blending alters the energy throughput capacity of pipeline systems.  This leads to significant changes in markets as well as practical operation of pipeline systems.  The electric power sector employs network optimization to compute location-specific prices for electricity based on the physics of energy flow \cite{litvinov2010design}, and similar market mechanisms for pricing in the pipeline transport of natural gas \cite{rudkevich2017hicss} and water networks \cite{biggar2022} were proposed.  Whereas the optimization-based economic analysis for existing electricity, water, and natural gas markets considered prices and quantities of commodities with homogeneous characteristics, such as heating value, a market analysis for a pipeline that carries blends of natural gas and hydrogen with multiple users would need to produce locational prices of natural gas and hydrogen for various suppliers, as well as prices of energy for downstream consumers that receive blends of various concentrations.  While optimization of gas pipeline flows with tracking of calorific values from different supply sources has been investigated \cite{hante2019complementarity}, the nonlinearities arising from hydrogen blending and need for practical market designs compel the new concepts developed here.

In this study, we use optimization subject to physical constraints that include flow equations, pressure limits, and compressor boost limits to evaluate hydrogen injection in a natural gas pipeline and the impact on the optimal flow allocation and pricing solution.  We extend recent results on the optimization formulations for pipeline flow allocation and capacity evaluation for gas mixtures \cite{sodwatana2023h2blend,zlotnik2023review}, and build upon the locational marginal pricing (LMP) concepts for energy networks \cite{rudkevich2017hicss} by examining heterogeneous gas flows in pipelines with incentives for offsetting carbon emissions by consuming hydrogen in end-use instead of natural gas. By solving for the dual variables, or the Lagrange multipliers, we derive the value of natural gas and hydrogen at each location in the network along with the decarbonization premium paid by end-users who consume gas that includes hydrogen.  These can be used to determine equitable subsidies or credits for hydrogen integration.  

The paper is organized as follows. In Section \ref{sec:pipe_modeling}, gas pipeline network modeling that includes heterogeneous physical flow  is introduced. Section \ref{sec:market_design} presents a market design concept and an optimization formulation that includes terms for carbon offset incentives. The key conceptual contribution of our study follows in Section \ref{sec:lagrange}, in which we examine the first-order optimality conditions and economic properties of market equilibria.  We demonstrate the application and scalability of the optimization model using 8-node and 40-node test networks in Section \ref{sec:case_study}, and show conditions for counter-intuitive market outcomes. 
We review the results of our case studies in Section \ref{sec:conclusion}.  Throughout the manuscript, we append SI units for variables in brackets after they are introduced.

\section{Heterogeneous Gas Flow Modeling} \label{sec:pipe_modeling}

We model a gas pipeline network as a connected and directed graph \((\mathcal{E},\mathcal{V})\), with physical junctions represented by nodes \(j \in \mathcal{V}\), and where each edge \((i,j) \in \mathcal{E}\) represents a pipe with flow from junction \(i\) to junction \(j\).  The subset \(\mathcal{C}\subset\mathcal{E}\) of edges contains compressors that boost gas pressure between pairs of nodes.  We also introduce the set \(\mathcal{G}\) of gNodes, following previously developed notation \cite{rudkevich2017hicss}. The gNodes correspond to individual participants in the pipeline market, which are either suppliers or consumers.   We use \(j(m) \in \mathcal{V}\) to indicate that \(j\) is the physical node location of a gNode \(m\). Each gNode \(m \in \mathcal{G}\) corresponds to a user of the pipeline network and can either be a hydrogen (H\(_2\)) supplier in the set \(\mathcal{G}^{H_2}_s\), a natural gas (NG) supplier in the set \(\mathcal{G}^{NG}_s\), or a consumer of gas in the set \(\mathcal{G}_d\).  For example, multiple gNodes may represent multiple suppliers or consumers at a physical node, or several gNodes can be used to represent various levels of market participation by a single entity at different price points. Finally, we specify a set of slack physical nodes \(\mathcal{V}_s\) that have nominal pressure value. Notations for sets in the pipeline network and indices are given in Table \ref{tab:sets_def}.

\begin{table}[h!]
\centering \normalsize
\begin{tabular}{|c|c|l|}
  \hline
  Index & Set  & Description \\ 
  \hline
  \hline
  \( j \) & \( \mathcal{V}\)  &   set of all physical nodes \\
  \( (i,j) \) & \( \mathcal{E}\)  & set of edges representing pipes\\
  \( (i,j) \) & \( \mathcal{C}\)  & set of edges representing compressors \\
  \( m \) & \( \mathcal{G}_s^{H_2}\) &  set of gNodes that inject hydrogen \\  
  \( m \) & \( \mathcal{G}_s^{NG}\)  & set of gNodes that inject natural gas \\
  \( m \) & \( \mathcal{G}_d\) & set of gNodes that withdraw gas \\
  \( j \) & \( \mathcal{V}_s\) &  set of slack physical nodes, subset of \(\mathcal{V}\) \\
  \hline
\end{tabular}
\caption{Notations for sets and indices. \label{tab:sets_def}}
\end{table}


Operational constraints are imposed on the nodal pressure \(P_j\), compressor boost ratio \(\alpha_{ij}\), the total mass flow \(\phi_{ij}\), and the mass fraction of hydrogen in junctions and along pipes, denoted by \(\gamma_j\) and \(\gamma_{ij}\) respectively. Supply and demand limits are imposed at the respective injection and withdrawal gNodes. 

\subsection{Pipe Flow Equations} \label{sec:pipe_equations}

We use the Weymouth equation to model the relation between mass flow and the pressures at the endpoints of a pipe \cite{rios2015optimization,kazi2024modeling}.  This relation is of the form 
\begin{equation}\label{eq:pipe_equations}
     P_{i}^{2} - P_{j}^{2} = \frac{f_{ij}L_{ij}}{D_{ij}A_{ij}^{2}} V_{ij} \phi_{ij} \left|\phi_{ij} \right|,  \quad \forall (i,j) \in \mathcal{E}, 
\end{equation}
where \(f_{ij}\), \(L_{ij}\), \(D_{ij}\) and \(A_{ij}\) are the friction factor, length, diameter, and cross-sectional area of pipe \((i,j)\), respectively. We suppose that the hydrogen concentration along the pipe is uniform and that flow is steady-state.  We use \(V_{ij}\) [(m/s)$^2$] to denote the square of the wave speed in the blended gas, approximated as a linear combination of the squared wave speeds \(a_{H_{2}}^{2}\) and \(a_{NG}^{2}\) [(m/s)$^2$] in hydrogen and natural gas as
\begin{equation} \label{eq:soundspeed}
    V_{ij} = \gamma_{ij}a_{H_{2}}^{2} + (1-\gamma_{ij})a_{NG}^{2}, \quad \forall (i,j) \in \mathcal{E}.
\end{equation}
In this study, we set \(a_{H_{2}}\) = 1090m/s and \(a_{NG}\) = 370m/s. Section \ref{sec:rescale} details the calculation for the wave speeds in hydrogen and natural gas.  The flow modeling in equations \eqref{eq:pipe_equations}-\eqref{eq:soundspeed} are the result of significant simplifications, as described in, e.g., Section V.A in \cite{roald2020uncertainty}. Equation \eqref{eq:pipe_equations} is the steady-state solution of the friction-dominated momentum conservation equation for pipeline flow.  Our focus in this study is on steady-state optimization of pipeline flow including hydrogen blending and the economic interpretation of the dual solution.  We leave precise quantitative calibration of the methodology, e.g., non-ideal gas modeling, to future studies.

\subsection{Nodal Compatibility Equations} \label{sec:node_equations}

The key nodal conditions for gas transport through a junction represent mass flow balance, which is linear for a homogeneous gas.  In the setting of blending multiple gases with concentration tracking, we require mass balance constraints at each physical node \(j\) for every gas constituent, such that  net incoming and outgoing flows of natural gas and hydrogen through each node are balanced. These mass balance equations depend on concentration and are thus nonlinear, of form
\begin{subequations} \label{eq:flowbalance}
\begin{equation} \label{eq:ngflowbalance}
    \!\!\! (1-\gamma_j) \sum_{k\in \partial_j^-} \phi_{jk} - \sum_{i\in \partial_j^+} (1-\gamma_{ij}) \phi_{ij} = \sum_{m\in \partial_j^g} s_m^{NG} -  (1-\gamma_j)\sum_{m\in \partial_j^g} d_m, \quad \forall j \in \mathcal{V}, 
\end{equation}
\begin{equation} \label{eq:h2flowbalance}
    \!\! \gamma_{j} \sum_{k\in\partial_{j}^{-}} \phi_{jk} - \sum_{i\in\partial_{j}^{+}} \gamma_{ij}\phi_{ij} = \! \sum_{m\in\partial_{j}^g} s_{m}^{H_2} - \gamma_{j} \! \sum_{m\in\partial_{j}^g} d_{m},  \quad \forall j \in \mathcal{V}, 
\end{equation}
\end{subequations}
where \(s_m^{NG}\) and \(s_m^{H_2}\) [kg/s] are the mass flow rate of natural gas and hydrogen at the injection gNode \(m\) of the physical node \(j\), and \(d_m\) [kg/s] is the mass flow rate of the delivered blended gas. Here, \(\partial_j^+\) and \(\partial_j^-\) are sets of nodes connected to \(j\) by incoming and outgoing edges, respectively, and $\partial_j^g$ denotes the set of gNodes associated with the physical node $j$. Adding equations \eqref{eq:ngflowbalance} and \eqref{eq:h2flowbalance} yields the total mass balance for the blended gas. To ensure continuity of gas concentration along the nodes and edges, we use a compatibility constraint 
\begin{equation}  \label{eq:continuity}
    \gamma_{ij} = \gamma_{i}, \quad \forall (i,j) \in \mathcal{E},
\end{equation}
which specifies that the hydrogen concentration of physical flow through edges \((i,j)\) must equal that at the sending node \(i\). Note that the converse is of course not true: input edge hydrogen concentrations can differ from each other and from the average nodal concentration due to blending at the node (where inputs are tracked).
We suppose that all input gas flows are mixed at nodes according to equations \eqref{eq:flowbalance} and flow as a homogeneous mixture along any output streams. The pressure \(\sigma_j\) at the slack nodes is
\begin{equation} \label{eq:slack_pressure}
    P_{j} = \sigma_j, \quad \forall j \in \mathcal{V}_s.
\end{equation}
At least one slack node is included following the convention in pipeline modeling for simulation to ensure a well-posed boundary value problem \cite{chaczykowski2018gas,brodskyi2024simulation}.  In formulating the nodal flow balance and mixing conditions in equations \eqref{eq:flowbalance} and \eqref{eq:continuity}, we assume that $\phi_{ij}\geq 0$.  The challenges in modeling for optimization of gas mixture flows on networks accounting for variable flow directions are described in our previous study \cite{kazi2024modeling}. Our assumption of fixed known flow in the oriented direction of graph edges avoids the need for mixed-integer constraints and streamlines our subsequent analysis of optimality conditions. This is a reasonable assumption in most operations currently, though possible future pipeline networks with meshed topologies may utilize more directional shifting, e.g., when hourly generation and consumption shifts significantly due to renewable-powered hydrogen injection at certain locations.

\subsection{Compressor Modeling} \label{sec:comp_modeling}

Compressor stations are used to maintain the flow of gas through transmission pipelines and compensate for pressure decrease in the direction of flow caused by friction. We apply a simplified model of such facilities, which may be complex sites with many machines, for our purpose of large-scale systems modeling.  We suppose that the action of a compressor station \((i,j) \in \mathcal{C}\) can be described by the boost ratio \(\alpha_{ij}\) between the suction and discharge pressures. The pressure at the end nodes \(i\) and \(j\) can then be related as
\begin{equation} \label{eq:comp_boost}
    P_{j} = \alpha_{ij} P_{i},   \quad \forall (i,j) \in \mathcal{C}.
\end{equation}
Higher hydrogen concentration in the gas blend requires more compression work to deliver the same amount of energy. To quantify the cost of compressor work, we first quantify the amount of work required to compress gas following traditional practice \cite{menon2005gas} using the adiabatic relation
\begin{equation} \label{eq:comp_power}
   \! W_c \!= \!\left(\! \frac{286.76\cdot (\kappa_{ij}-1)\cdot T}{G_{ij}{\kappa_{ij}}} \!\right) \!\left( \alpha_{ij}^{m} \!-\! 1 \right) \left|\phi_{ij}\right|, \quad \forall (i,j) \!\in\! \mathcal{C},
\end{equation}
where \(T\) [K] is the temperature of gas entering the compressor, and is nominally 288.7 K in our study. Here, \(\kappa_{ij}\) and \(G_{ij}\) denote the specific heat capacity ratio and the specific gravity of the gas mixture flowing through the compressor, and \(m = (\kappa_{ij} - 1)/\kappa_{ij}\).  These parameters are $\kappa_{NG}=1.304$ and $G_{NG}=0.5537$ for methane and $\kappa_{NG}=1.405$ and $G_{NG}=0.0696$ for hydrogen.  
When using the equation \eqref{eq:comp_power} to approximate compressor work,
we set \(\kappa\) and \(G\) to nominal constant values \(\kappa_{nom} = 1.308\) and \(G_{nom} = 0.574\) that correspond to a hydrogen mass fraction of $\gamma_{ij}\approx 0.05$ to avoid excessive complexity in the compressor cost function.  Using the above expression, we assess the economic cost of actuating gas flow through the pipeline by using \(\eta\) [\$/kw-s] to represent the cost of operating an electric drive gas compressor. In our examples, we use \(\eta\) = \$0.13/3600kw-s, which is equivalent to an electricity price of 13$\mbox{\textcent}$/kWh. This results in
\begin{equation} \label{eq:compsimplified}
    \eta W_c = \eta K \left( \alpha_{ij}^{m_{nom}}-1 \right) \left|\phi_{ij}\right|, \quad \forall (i,j) \in \mathcal{C},
\end{equation}
where $K$ denotes a constant that represents the quantity in parentheses in equation \eqref{eq:comp_power}. 
The above simplifications result in values of \(K\) = 22.18 and \(m\) = 0.325.  The parameters  \(\kappa\) and \(G\) are set to nominal values in order to simplify the optimization, and more practical justification for this is given in Section \ref{sec:optim_problem}.

\subsection{Pressure, Compressor, and Concentration Limits} \label{sec:engineering_limits}

In addition to the equality constraints in equations \eqref{eq:pipe_equations}-\eqref{eq:compsimplified}, inequality constraints are imposed on the pressure, compressor boost ratio, and hydrogen concentration to reflect the engineering, operating, and contractual limitations that the gas transmission pipeline is subject to. We suppose that the minimum operational pressure limit and the minimum and maximum hydrogen concentration are specified at each node:
\begin{subequations} \label{eq:nodallimits}
\begin{align} 
    P_{j}^{min} \le P_{j}  \quad & \forall j \in \mathcal{V}, \label{eq:minpressure} \\
    \gamma_{j}^{min} \le \gamma_{j} \le \gamma_{j}^{max} \quad & \forall j \in \mathcal{V}. \label{eq:conclimits}
\end{align}
\end{subequations}
We suppose that compressors can only increase pressure because regulation to reduce pressure is typically not done along midstream pipelines but at citygates to local distribution systems. Therefore, we suppose that \(\alpha_{ij}\) has a lower bound of 1. In addition, the compressor boost ratio is bounded by the maximum allowable pressure \(P_{j}\) at the discharge node and the maximum boost ratio \(\alpha_{ij}^{max}\):
\begin{subequations} \label{eq:complimits}
\begin{align} 
    \alpha_{ij}P_{i} \le P_{j}^{max}  \quad & \forall (i,j) \in \mathcal{C}, \label{eq:maxpressure} \\
    1 \le \alpha_{ij} \le \alpha^{max}_{ij} \quad & \forall (i,j) \in \mathcal{C}. \label{eq:compratiolimits}
\end{align}
\end{subequations}
The above constraints define the physical state of the pipeline system.  The next section includes additional constraints and an objective function that will be used to define an optimization problem to clear a double-sided single auction market.

\section{Optimization Formulation for Market Design with Incentives} \label{sec:market_design}

We design an auction market mechanism for natural gas pipeline capacity following inspiration from an early study published by the U.S. Federal Energy Regulatory Commission \cite{united1987gas}, which was revisited in a more recent study on the concept of locational pricing for natural gas transmission \cite{rudkevich2017hicss}.  Here, we include additional details to account for hydrogen blending in the capacity market, as well as incentives for the carbon dioxide emissions reduction that results from the displacement of end-use natural gas by green hydrogen.  The mechanism is designed such that the objective can represent the economic value generated by the pipeline system for its users, who provide offers to sell and bids to buy energy in the form of the delivered gas mixture. The objective includes the cost to transport gas using gas compressors as discussed in Section \ref{sec:comp_modeling}, incentives for reducing carbon emissions by using hydrogen to displace natural gas combustion, and the revenue collected because of differences between buyer and seller prices for the gas mixture arising due to network congestion.

\subsection{Supply and Demand Limits} \label{sec:economic_limits}

To formulate the proposed market structure, we suppose that each seller of hydrogen gives a price and quantity offer (\(c_m^{H_2}\), \(s_m^{max,H_2}\)), with units ([\$/kg], [kg/s]), at the corresponding injection gNode. Similarly, each natural gas seller provides a price and quantity offer (\(c_m^{NG}\), \(s_m^{max,NG}\)), with units ([\$/kg], [kg/s]). The quantity components of offers from suppliers are used to formulate the supplier-side constraints 
\begin{subequations} \label{eq:supplylimits}
\begin{align}
    0\le s_{m}^{NG} \le s_{m}^{max,NG}, \quad & \forall m \in \mathcal{G}_s^{NG}, \label{eq:supplylimits_ng} \\
    0\le s_{m}^{H_2} \le s_{m}^{max,H_2}, \quad & \forall m \in \mathcal{G}_s^{H_2},  \label{eq:supplylimits_h2}
\end{align}
\end{subequations}
where the quantity is used as the upper bound value on the supply, and \(s_m^{H_2}\) and \(s_m^{NG}\) denote the optimized injection flows for hydrogen and natural gas, respectively. We suppose that each flexible customer places a bid consisting of the price and quantity (\(c_m^{d}\), \(g_m^{max}\)) for energy, with units ([\$/MJ], [MJ/s]). There may also be customers with fixed demand \(\bar{g}_m\) [MJ/s] that is pre-determined outside the market mechanism, whose energy consumption must be served by the pipeline regardless of the market outcome. The demand-side constraints are
\begin{subequations}\label{eq:energydemand}
\begin{align}
    0\le d_{m} \left( R_{H_{2}}  \gamma_{j(m)} + R_{NG} (1\!-\gamma_{j(m)}) \right) \!\le\! g_{m}^{max}, & \quad   \forall m \!\in \mathcal{G}_{d,o}, \label{eq:demand_opt} \\
   d_{m} \left( R_{H_{2}}  \gamma_{j(m)} + R_{NG} (1-\gamma_{j(m)}) \right) = \bar{g}_{m}, & \quad   \forall m \in  \mathcal{G}_{d,f}, \label{eq:demand_fixed}
\end{align}
\end{subequations}
where \(d_m\) [kg/s] is the optimized withdrawal mass flow rate of the blended gas. The demands of flexible customers at gNodes in the set $\mathcal{G}_{d,o}$ are limited by constraint \eqref{eq:demand_opt} while consumption by customers at gNodes $\mathcal{G}_{d,f}$ with fixed demand are subject to constraint \eqref{eq:demand_fixed}.  Because the price and quantity bids are in units of equivalent power, we use the expression
\begin{equation}\label{eq:energy_conversion}
  R(\gamma_{j(m)}) =  R_{H_{2}}  \gamma_{j(m)} + R_{NG} (1\!-\gamma_{j(m)}),
\end{equation}
which gives the calorific value of the gas based on the hydrogen mass fraction, where \(R_{H_2}\) = 141.8 MJ/kg and \(R_{NG}\) = 44.2 MJ/kg are the calorific values for hydrogen and natural gas. The quantity $R(\gamma_{j(m)})$ gives the composition-dependent conversion factor between mass flow and energy flow. We will henceforth be using the shorthand notation \(R(\gamma_{j(m)})\).

\subsection{Carbon Offset Incentives} \label{sec:emissions_modeling}

A straightforward method for creating incentives to lower carbon dioxide (CO\(_2\)) emissions in optimization-based markets is to add a term to the objective function that quantifies the value of avoided emissions.  We suppose that \(c_m^{CO_2}\) [\$/kgCO\(_2\)] is an incentive paid to the market administrator for CO\(_2\) emissions that are avoided when consumers use hydrogen instead of natural gas to produce the same unit of energy.  These avoided carbon emissions are denoted by \(E_m\) [kgCO\(_2\)/s], and are approximated as
\begin{equation} \label{eq:carbon_offset}
    E_m = d_m \gamma_{j(m)} \cdot \frac{R_{H_2}}{R_{NG}} \cdot \zeta, \quad  \forall m \in \mathcal{G}_d, 
\end{equation}
where \(\zeta\) = 44/16 is the approximate ratio of the molecular weight of carbon dioxide to methane.  The incentives received by the market administrator can be passed through to consumers who pay more for energy because of hydrogen blending, which we refer to as the locational decarbonization premium.  We will examine this concept in Section \ref{subsec:economic}. 

\subsection{Objective Function and Optimization Formulation} \label{sec:optim_problem}

The objective function for the proposed optimization-based market mechanism is to maximize the economic value produced by operating the pipeline system.  The three components of the objective function quantify (i) the market revenue collected by selling energy to consumers minus the cost of procuring natural gas and hydrogen from suppliers; (ii) the carbon emissions mitigation incentives collected by the market operator on behalf of all participants; and (iii) the cost of operating gas compressors.  We suppose that these objective function components are, respectively, of form
\begin{subequations} \label{eq:obj_components}
\begin{align}
    J_{MR} & = \sum_{m\in G} (c_m^{d} d_m R(\gamma_{j(m)}) -   c_m^{H_2} s_m^{H_2}  - c_m^{NGs} s_m^{NG}), \label{eq:obj_components_MR}\\
    J_{CEM} & =  \sum_{m\in G} c_m^{CO_2} E_m, \label{eq:obj_components_CEM}\\
    J_{GC} & = \eta \sum_{c\in C} W_c =  \eta \sum_{c\in C} K \left( \alpha_{ij}^{m_{nom}}-1 \right) \left|\phi_{ij}\right|, \label{eq:obj_components_GC}
\end{align}
\end{subequations}
and combining these components to maximize the market revenue and carbon emissions mitigation while minimizing compressor energy yields
\begin{equation}\label{eq:obj}
\begin{split}
    J_{EV} & = J_{MR} + J_{CEM} - J_{GC} \\ & = \sum_{m\in G} \bigg(c_m^{d} d_m R(\gamma_{j(m)}) -   c_m^{H_2} s_m^{H_2}  - c_m^{NGs} s_m^{NG} +  c_m^{CO_2} E_m \bigg)  - \eta \sum_{c\in C} W_c.
\end{split}
\end{equation}
We note that the energy used to move gas through a pipeline system, represented here as $J_{GC}$, would be two or three orders of magnitude smaller than the economic value produced by transporting energy from suppliers to consumers, as quantified by $J_{MR}$.  The approximate modeling of adiabatic gas compression described in Section \ref{sec:comp_modeling} can be justified because in practice the value of $J_{GC}$ will  not significantly affect the optimal solution, while it does provide regularization to prevent degeneracy, e.g., a set of of feasible points that yields the same value of $J_{MR}$. The latter situation could arise when no hydrogen is injected, in which case equation \eqref{eq:obj_components_MR} would be linear.

We formulate the optimization problem by maximizing the objective \eqref{eq:obj} subject to the above constraints that define network flows, engineering limits, and market conditions.  The resulting optimization formulation is
\begin{equation} \label{eq:optimization_formulation}
\begin{array}{ll}
    \!\!\!\! \mathrm{max}  & J_{EV} \triangleq \text{max economic value objective } \eqref{eq:obj} \\
    \!\!\!\! \text{s.t.} & \text{pressure balance } \eqref{eq:pipe_equations} \\
    & \text{NG flow balance } \eqref{eq:ngflowbalance} \\
    & \text{H\textsubscript{2} flow balance } \eqref{eq:h2flowbalance} \\
    & \text{H\textsubscript{2} concentration constraints } \eqref{eq:continuity},\eqref{eq:conclimits} \\
    & \text{pressure constraints } \eqref{eq:slack_pressure},  \eqref{eq:minpressure} \\
    & \text{compressor boost limits } \eqref{eq:complimits} \\
    & \text{supply limits } \eqref{eq:supplylimits} \\
    & \text{demand limits } \eqref{eq:energydemand}
\end{array}
\end{equation}
with all decision variables modeled as real and continuous. A comprehensive nomenclature for the decision variables, derived quantities, physical parameters, and operational parameters used to formulate the above problem is provided in \ref{app:nomenclature}. 

Observe that in principle, the carbon emissions mitigation incentive \(c_m^{CO_2}\) could be made node-dependent, and this could be imposed by an external entity, determined by the market administrator, or provided as part of a consumer bid.  For example, the locational value of the incentive could vary because of differences in  incentives in regions through which a pipeline passes.   While we index the incentive \(c_m^{CO_2}\) by the gNode $m\in\mathcal{G}$, we set it to a uniform network-wide value for each of the  scenarios that we examine.  Similarly, in this setting we suppose that customers are not in general able to calibrate the mass fraction of hydrogen in the gas delivered to their location, but rather the mixture fraction is governed by the physics of mixing, the location relative to the injection point, and the economics of the resulting formulation driven by hydrogen blending. However, locational blending requirements or limitations at a node $j\in\mathcal{V}$ could be specified by adjusting $\gamma_{j}^{min}$ and/or $\gamma_{j}^{max}$, although setting $\gamma_{j}^{min}>0$ may preclude the existence of feasible solutions to the mathematical program defined below.  In the scenarios examined here, we consider uniform network-wide values of $\gamma_{j}^{max}$ and $\gamma_{j}^{min}\equiv 0$.  We leave analysis of more complex, non-uniform locational objectives for decarbonization and locational constraints on hydrogen mass fraction for future studies. In the following section, we examine first-order optimality conditions for problem \eqref{eq:optimization_formulation}.

\section{Optimality Conditions and Economic Analysis}  \label{sec:lagrange}

Markets for the procurement and transport of natural gas through large-scale pipeline systems are typically managed using frameworks that are not designed for responsive price formation \cite{macavoy2008natural}.   There is a physics-based optimization-driven market that uses a linearized gas flow model that has been used for a regional network in Australia \cite{pepper2012implementation}, but this is an exception.  Holders of firm transportation rights typically sell unused capacity, which is bundled with gas supply and traded bilaterally in a daily spot market.  In contrast, it is standard to use optimization to form forward and real-time prices in wholesale electricity markets \cite{litvinov2010}, and prices posted by system operators drive participation of generation managers \cite{manual2012energy}.  In particular, the concept of locational marginal prices has been widely implemented in practice for power systems \cite{schweppe2013spot}. Here we examine mathematical programming as a mechanism for responsive, optimal scheduling and price formation that can also facilitate the incorporation of new paradigms for pipeline systems such as hydrogen blending and decarbonization initiatives.  Specifically, we consider optimality conditions and dual solutions for problem \eqref{eq:optimization_formulation}, and examine the potentials and challenges for its use in price formation. 

Following a previously proposed optimization-based formulation for a double auction market for natural gas pipeline capacity \cite{rudkevich2017hicss}, we examine the first-order optimality conditions for problem \eqref{eq:optimization_formulation}. First, we derive the Lagrangian of the optimization problem \eqref{eq:optimization_formulation} by adjoining the constraints as listed in the problem using Lagrange multipliers to the objective function expressed as a minimization.  This results in 

\begin{align}\label{eq:full_lagrangian}
    \mathcal{L} &= \!\sum_{m\in G} \bigg(\!\!-c_m^{d} d_m R(\gamma_{j(m)}) \!+\!  c_m^{H_2} s_m^{H_2}  \!+\! c_m^{NGs} s_m^{NG} \!-\!  c_m^{CO_2} d_m \gamma_{j(m)} \frac{R_{H_2}}{R_{NG}}   \zeta \bigg)   \!+\! \sum_{i,j\in C} K \left( \alpha_{ij}^{m_{nom}}\!-\!1 \right) \left|\phi_{ij}\right| \nonumber   \\
    & \quad + \!\sum_{j\in \mathcal{V}} \lambda_j^{NG} \bigg((1-\gamma_j) \textstyle\sum_{k\in \partial_j^-} \phi_{jk} - \sum_{i\in \partial_j^+} (1-\gamma_{ij}) \phi_{ij} - \sum_{m\in \partial_j^g} s_m^{NG} \!+\! (1-\gamma_j)\sum_{m\in \partial_j^g} d_m \!\bigg) \nonumber  \\
    & \quad + \!\sum_{j\in \mathcal{V}} \lambda_j^{H_2} \bigg( \gamma_{j} \textstyle\sum_{k\in\partial_{j}^{-}} \phi_{jk} - \sum_{i\in\partial_{j}^{+}} \gamma_{ij}\phi_{ij} - \sum_{m\in\partial_{j}^g} s_{m}^{H_2} + \gamma_{j}  \sum_{m\in\partial_{j}^g} d_{m}  \bigg) \nonumber  \\
    & \quad + \!\sum_{(i,j)\in \mathcal{E}} \mu_{ij} \bigg( P_i^2 - P_j^2 - \frac{f_{ij}L_{ij}}{D_{ij}A_{ij}^2} \big(\gamma_{ij} a_{H_2}^2 + (1 - \gamma_{ij}) a_{NG}^2\big) \phi_{ij} |\phi_{ij}| \bigg) \nonumber  \\
    & \quad  + \!\sum_{(i,j)\in \mathcal{E}} \omega_{ij}^e \left( \gamma_j - \gamma_{ij}\right) + \sum_{j\in \mathcal{V}} \left( \omega_j^l \left( \gamma_j^{\min} - \gamma_{j} \right) + \omega_j^u \left( \gamma_j - \gamma_j^{\max}\right) \right)  + \sum_{j\in \mathcal{V}} \beta_j^l \left( P_j^{\min} - P_j \right) \nonumber \\
    & \quad + \!\sum_{j\in \mathcal{V}_s} \beta_j^e \left( P_j - \sigma_j \right) 
     + \sum_{(i,j)\in \mathcal{C}} \theta_{ij}^e \left( P_j^2 - \alpha_{ij}^2 P_i^2 \right) + \sum_{(i,j)\in \mathcal{C}} \theta_{ij}^u \left( \alpha_{ij} P_i - P_{j}^{\max} \right)   \nonumber\\
     & \quad + \!\sum_{(i,j)\in \mathcal{C}}\!\!\left( \theta_{ij}^{c,l} \left( 1-\alpha_{ij} \right) \!+\! \theta_{ij}^{c,u} \left( \alpha_{ij} - \alpha_{ij}^{\max} \right) \right) 
     + \!\!\!\sum_{m\in \mathcal{G}_s^{NG}} \!\!\!\left( \chi_m^{NG,l} \left(  - s_m^{NG} \right) \!+\! \chi_m^{NG,u} \left( s_m^{NG} - s_{m}^{\max,NG} \right) \right)  \quad \nonumber \\
     & \quad + \sum_{m\in \mathcal{G}_s^{H_2}} \!\left( \chi_m^{H_2,l} \left(  - s_m^{H_2} \right) + \chi_m^{H_2,u} \left( s_m^{H_2} - s_{m}^{\max,H_2} \right) \right)    
     + \sum_{m\in \mathcal{G}_{d,o}} \chi_m^l \left( -d_m R(\gamma_{j(m)})  \right) \nonumber\\
     & \quad + \sum_{m\in \mathcal{G}_{d,o}} \chi_m^u \left( d_m R(\gamma_{j(m)}) - g_m^{\max} \right) \quad + \sum_{m\in \mathcal{G}_{d,f}} \chi_m^f \left( d_m R(\gamma_{j(m)}) - \bar{g}_m \right). 
\end{align}
In the Lagrangian in equation \eqref{eq:full_lagrangian}, we use \(\lambda_j^{NG}\) and \(\lambda_j^{H_2}\) to denote the Lagrange multipliers corresponding to the flow balance constraints \eqref{eq:ngflowbalance} and \eqref{eq:h2flowbalance}, respectively, and \(\mu_{ij}\) is the Lagrange multiplier of the pressure balance constraint \eqref{eq:pipe_equations}. The multipliers corresponding to the pressure limits are \(\beta_j^{l}\) and \(\beta_j^{e}\), while the multipliers for the compressor boost limits are \(\theta_{ij}^{e}, \theta_{ij}^{u}, \theta_{ij}^{c,l}\) and \(\theta_{ij}^{c,u}\). Here \(\omega_j^{l}, \omega_j^{u}\) and \(\omega_{ij}^{e}\) are the Lagrange multipliers corresponding to the hydrogen concentration limits \eqref{eq:conclimits} and the continuity constraint \eqref{eq:continuity}. \(\chi_m^{NG,l}, \chi_m^{NG,u}, \chi_m^{H_2,l}\) and \(\chi_m^{H_2,u}\) are the multipliers that correspond to the lower and upper bounds of the gas supply constraints \eqref{eq:supplylimits}, while \(\chi_m^{d}\) and \(\chi_m^{g}\) are the multipliers of the withdrawal constraints \eqref{eq:energydemand}. All the multipliers that are used to synthesize the Lagrangian in equation \eqref{eq:full_lagrangian} are tabulated in Table \ref{tab:multipliers_defs}, together with the corresponding constraint and dimensional units. Certain Lagrange multipliers can be interpreted as marginal prices, which can be decomposed into price components of interest, as will be discussed subsequently.

\subsection{Optimality Conditions} \label{subsec:optimality}

We now evaluate the Karush-Kuhn-Tucker (KKT) conditions that optimal solutions to the nonlinear program \eqref{eq:optimization_formulation} must satisfy \cite{peressini1988mathematics,bertsekas2014constrained}. Specifically, we examine first-order conditions for stationarity, primal feasibility, dual feasibility, and complementary slackness, as well as some of the second-order sufficiency conditions.

\begin{table}[t!]
\centering \normalsize
\begin{tabular}{|c|c|c|}
  \hline
  Multiplier & Constraint  & Units \\ 
  \hline
  \hline
  $\lambda_j^{NG}$  & NG flow balance \eqref{eq:ngflowbalance}  &  \$/kg NG \\
  $\lambda_j^{H_2}$ & H$_2$ flow balance \eqref{eq:ngflowbalance} & \$/kg H$_2$\\
  $\mu_{ij}$ & pipe equation \eqref{eq:pipe_equations}  & $\$/ (s\cdot$Pa$^2$) \\
  $\omega_{ij}^{e}$ & concentration continuity \eqref{eq:continuity} & \$/s \\  
  $\omega_j^{l},\, \omega_j^{u}$ & concentration limits \eqref{eq:conclimits} & \$/s \\
  $\beta_j^{l}$ & minimum pressure \eqref{eq:minpressure} & $\$/ (s\cdot$Pa) \\
  $\beta_j^{e}$ & slack node pressure \eqref{eq:slack_pressure} & $\$/ (s\cdot$Pa) \\
  $\theta_{ij}^{e}$ & compressor boost \eqref{eq:comp_boost}& $\$/ (s\cdot$Pa$^2$) \\
  $\theta_{ij}^{u}$ & discharge pressure \eqref{eq:maxpressure} & $\$/ (s\cdot$Pa) \\
  $\theta_{ij}^{c,l}, \, \theta_{ij}^{c,u}$ & compressor ratio limits \eqref{eq:compratiolimits} & \$/s \\
  $\chi_m^{NG,l},\, \chi_m^{NG,u}$ & NG supplier limits \eqref{eq:supplylimits_ng} & \$/kg NG \\
  $\chi_m^{H_2,l}\,\chi_m^{H_2,u}$ & H$_2$ supplier limits \eqref{eq:supplylimits_h2} & \$/kg H$_2$ \\
  $\chi_m^{l},\, \chi_m^{u}$ & optimized demand \eqref{eq:demand_opt} & \$/MJ \\
  $\chi_m^{f}$ & fixed demand \eqref{eq:demand_fixed} & \$/MJ \\
  \hline
\end{tabular}
\caption{Lagrange multipliers, corresponding constraints, and units. \label{tab:multipliers_defs}}
\end{table}

The complementary slackness conditions for the inequality constraints are listed below in equations \eqref{eq:csc}, together with the indexing in the network component sets.
\begin{subequations} \label{eq:csc} 
\begin{alignat}{2}
\quad & \beta_j^e \left( P_j^{\min} - P_j \right) = 0 \quad & \forall j \in \mathcal{V} \label{eq:csc1}\\
\quad & \theta_{ij}^u \left( \alpha_{ij} P_i - P_{ij}^{\max} \right) = 0 \quad & \forall (i,j) \in \mathcal{E} \label{eq:csc2} \\
\quad & \theta_{ij}^{c,l} \left( 1-\alpha_{ij} \right) = 0 \quad & \forall (i,j) \in \mathcal{E} \label{eq:csc3} \\
\quad & \theta_{ij}^{c,u} \left( \alpha_{ij} - \alpha_{ij}^{\max} \right) = 0 \quad & \forall (i,j) \in \mathcal{E} \label{eq:csc4} \\
\quad & \omega_j^l \left( \gamma_j^{\min} - \gamma_{j} \right) = 0 \quad & \forall j \in \mathcal{V} \label{eq:csc5}\\
\quad & \omega_j^u \left( \gamma_j - \gamma_j^{\max}\right) = 0 \quad & \forall j \in \mathcal{V} \label{eq:csc6} \\
\quad & \chi_m^{H_2,l} \left(  - s_m^{H_2} \right) = 0 \quad & \forall m \in \mathcal{G}_s^{H_2} \label{eq:csc7} \\
\quad & \chi_m^{H_2,u} \left( s_m^{H_2} - s_{m}^{\max,H_2} \right) = 0 \quad & \forall m \in \mathcal{G}_s^{H_2} \label{eq:csc8} \\
\quad & \chi_m^{NG,l} \left( - s_m^{NG} \right) = 0 \quad & \forall m \in \mathcal{G}_s^{NG} \label{eq:csc9} \\
\quad & \chi_m^{NG,u} \left( s_m^{NG} - s_{m}^{\max,NG} \right) = 0 \quad & \forall m \in \mathcal{G}_s^{NG} \label{eq:csc10} \\
\quad & \chi_m^l \left( - d_m R(\gamma_{j(m)})\right) = 0 \quad & \forall m \in \mathcal{G}_{d,o} \label{eq:csc11} \\
\quad & \chi_m^u \left( d_m R(\gamma_{j(m)}) - g_m^{\max}\right)  = 0 \quad & \forall m \in \mathcal{G}_{d,o} \label{eq:csc12}
\end{alignat}
\end{subequations}
The multiplier sign condition requires positive multipliers for inequality constraints, i.e.,
\begin{subequations} \label{eq:msc}
\begin{alignat}{2}
\beta_j^e,\, \theta_{ij}^u,\, \theta_{ij}^{c,l},\, \theta_{ij}^{c,u},\, \omega_j^l,\, \omega_j^u \geq 0, \\
\chi_m^{H_2,l},\, \chi_m^{H_2,u},\, \chi_m^{NG,l},\, \chi_m^{NG,u},\, \chi_m^l, \chi_m^u  \ge 0. \label{eq:msc1}
\end{alignat}
\end{subequations}
We can then take the first derivative of the Lagrangian \(\mathcal{L}\) in equation \eqref{eq:full_lagrangian} with respect to each primal variable, to obtain the first-order necessary conditions for optimality. The conditions are listed in equations \eqref{eq:ldc} below, preceded by the primal variable with respect to the variation evaluated.  The first order derivative conditions are

\begin{subequations} \label{eq:ldc}
\begin{alignat}{2}
s_m^{NG}: \quad & 0 = c_m^{NGs} - \lambda_{j(m)}^{NG} - \chi_m^{NGl} + \chi_m^{NGu}  \quad & \forall m \in \mathcal{G}_s^{NG} \label{eq:ldc1}\\
s_m^{H_2}: \quad & 0 = c_m^{H_2} - \lambda_{j(m)}^{H_2} - \chi_m^{H_2l} + \chi_m^{H_2u}  \quad & \forall m \in \mathcal{G}_s^{H_2} \label{eq:ldc2}\\
d_m: \quad & 0 = -c_m^{d} R(\gamma_{j(m)}) - c_m^{CO_2}  \gamma_{j(m)}  \frac{R_{H_2}}{R_{NG}}   \zeta   + \lambda_{j(m)}^{NG} (1-\gamma_j) + \lambda_{j(m)}^{H_2} \gamma_j  \nonumber\\ & \qquad - \chi_m^lR(\gamma_{j(m)})   + \chi_m^uR(\gamma_{j(m)}) + \chi_m^fR(\gamma_{j(m)}) \quad & \forall m \in \mathcal{G}_d \label{eq:ldc3} \\
\phi_{ij}: \quad & 0 = (\lambda_i^{NG}(1-\gamma_i) - \lambda_j^{NG}(1-\gamma_{ij})) + (\lambda_i^{H_2}\gamma_{i} - \lambda_j^{H_2}\gamma_{ij})\nonumber \\  & \qquad - 2\mu_{ij} \frac{f_{ij}L_{ij}}{D_{ij}A_{ij}^2} (\gamma_{ij}a_{H_2}^2 + (1-\gamma_{ij})a_{NG}^2)\left|\phi_{ij}\right| \quad & \forall (i,j) \in \mathcal{E} \label{eq:ldc4} \\
\phi_{ij}: \quad & 0 = (\lambda_i^{NG}(1-\gamma_i) - \lambda_j^{NG}(1-\gamma_{ij})) + (\lambda_i^{H_2}\gamma_{i} - \lambda_j^{H_2}\gamma_{ij}) \nonumber \\ \quad & \qquad + K \left( \alpha_{ij}^{m_{nom}}-1 \right)  \quad & \forall (i,j) \in \mathcal{C} \label{eq:ldc5} \\
\alpha_{ij}: \quad &  0 = K |\phi_{ij}| m_{nom} \alpha_{ij}^{m_{nom}-1} -2\theta_{ij}^e\alpha_{ij}P_i^2 + \theta_{ij}^u P_i - \theta_{ij}^{c,l} + \theta_{ij}^{c,u} \quad & \forall (i,j) \in \mathcal{C} \label{eq:ldc6} \\
\gamma_{ij}: \quad & 0 = (\lambda_j^{NG} - \lambda_j^{H_2}) \phi_{ij} -\mu_{ij} \frac{f_{ij}L_{ij}}{D_{ij}A_{ij}^2} \phi_{ij} \left|\phi_{ij}\right| (a_{H_2}^2 - a_{NG}^2) - \omega_{ij}^e \quad & \forall (i,j) \in \mathcal{E}  \label{eq:ldc7} \\
\gamma_{j}: \quad & 0 = \!\sum_{m\in \partial_j^g} \bigg(-c_m^{d}d_m (R_{H_2} - R_{NG}) +  c_m^{CO_2} d_m \frac{R_{H_2}}{R_{NG}}   \zeta  - \lambda_j^{NG}d_m \bigg)  \nonumber \\ \quad & \quad + \sum_{m\in \partial_j^g} \!\left(\lambda_j^{H_2} d_m  \!+\! (\chi_m^u \!+\! \chi_m^f) d_m( R_{H_2} \!-\! R_{NG}) \!-\! \chi_m^l d_m( R_{H_2} \!-\! R_{NG}) \right) \nonumber \\ & \qquad - \lambda_j^{NG}\phi_{jk}   +  \lambda_j^{H_2}\phi_{jk} - \omega_j^l + \omega_j^u + \omega_{ij}^e \quad & \forall j \in \mathcal{V} \label{eq:ldc8} \\
P_{j}: \quad & 0 = 2\mu_{jk}P_j - 2\mu_{ij}P_j - \beta_j^l + \beta_j^e + 2\theta_{ij}^e P_j - 2\theta_{jk}^e \alpha_{jk}^2 P_j + \theta_{jk}^u\alpha_{jk} \quad & \forall j \in \mathcal{V} \label{eq:ldc9}
\end{alignat}
\end{subequations}
Finally, we consider the second-order sufficient conditions (SOSC) evaluated for the problem \eqref{eq:optimization_formulation} with respect to the Lagrangian \eqref{eq:full_lagrangian}. The SOSC are \cite{bertsekas2014constrained}:
\begin{subequations}
    \begin{align}
    d^T \nabla^2_{xx} \mathcal{L} (x^*, \lambda^*) d \geq 0, \qquad \forall d \in \mathcal{T}_F(x^*),
    \end{align}
    where $x^*$ and $\lambda^*$ denote the primal and dual solution, respectively. The direction vector $d$ needs to be in the tangent cone $\mathcal{T}_F(x^*)$ of the feasible region, which is defined as
    \begin{align}
        \mathcal{T}_F(x^*)  = & \{ \,\, d \,\, | \,\, 
        \nabla g_i(x^*)^T d \leq 0 \,\, \forall i \in I_g(x^*), \lambda^*_i = 0, \\ 
        & \,\, \qquad \nabla g_i(x^*)^T d = 0 \,\, \forall i \in I_g(x^*), \lambda^*_i > 0, \\ & \,\, \qquad  
        \nabla h(x^*)^T d = 0 \},
    \end{align}
    where $g$ and $h$ denote the inequality and equality constraints respectively, and $I_g$ is the active set of inequalities for the point $x^*$.
\end{subequations}

We focus on the specific variables $\phi_{ij}$ when examining the SOSC. Observe that none of the inequality constraints \eqref{eq:csc} include $\phi_{ij}$, and the only equality constraints involving $\phi_{ij}$ are the pressure balance equation \eqref{eq:pipe_equations} and individual mass balance equations \eqref{eq:flowbalance}. Thus, the direction vector element $d_{\phi}$ corresponding to the variable $\phi_{ij}$ is related to the elements of $d$ corresponding to $\gamma_{ij}$, $P_i$, $P_j$, and $s^{H_2}_m$, $s^{NG}_m$, and $d_m$ through the gradients $\nabla h$ where $h$ denotes the equality constraints in equations \eqref{eq:pipe_equations} and \eqref{eq:flowbalance}. For the $d_{\phi} = 1$, we can solve for the other elements of the direction vector and verify the SOSC. Because the maximum degree of any variable in the Lagrangian $\mathcal{L}$ is 2 for variable $\phi_{ij}$, related to equation \eqref{eq:pipe_equations}, the hessian of the Lagrangian $\nabla^2_{xx} \mathcal{L}$ has only one non-zero element on its diagonal. Thus the SOSC for the non-convex problem can be formulated as
\begin{equation} \label{eq:second_order_flow}
 - 2\mu_{ij} \frac{f_{ij}L_{ij}}{D_{ij}A_{ij}^2} (\gamma_{ij}a_{H_2}^2 + (1-\gamma_{ij})a_{NG}^2)\cdot\mathrm{sign}(\phi_{ij}) > 0, \quad \forall (i,j) \in \mathcal{E}.
\end{equation}
Based on our assumption that $\phi_{ij}\geq0$, we may conclude that $\mu_{ij}\leq0$. There are a collection of additional regularity conditions, or constraint qualifications, that can be used to determine whether a given stationary point that satisfies the stationarity conditions \eqref{eq:ldc}.  These include linear independence constraint qualification (LICQ) that active inequality constraints and gradients of the equality constraints must be linearly independent, and Slater's condition for convex problems. For general nonlinear constrained optimization, including conic or quadratic problems, there may be a duality gap \cite{udell2016bounding}. Rigorous analysis and verification of the duality gap, or dual Lagrange problem, for solutions to formulation \eqref{eq:optimization_formulation} is beyond our present scope.  Numerous studies on duality theory for extremal problems \cite{ioffe1968duality,rockafellar2009variational} and analysis of Lagrange multipliers in constrained nonlinear optimization \cite{borwein2016variational,bertsekas2014constrained} describe the challenges in this setting. 
Our subsequent economic interpretations of primal-dual solutions to problem \eqref{eq:optimization_formulation} and the associated KKT conditions assume a zero or negligible duality gap. We also verify the derived relations using computational case studies.

\subsection{Economic Interpretations} \label{subsec:economic}

The concepts of spot pricing and congestion rents are well-developed in electricity markets \cite{schweppe2013spot}, and similar concepts have been developed for gas pipelines \cite{rudkevich2017hicss}.  A critical issue in the design of markets, particularly involving transportation and/or networks, is revenue adequacy.  This notion was first suggested in the context of railroad markets, where revenue adequacy refers to a higher return on investment compared to the cost of capital.   There has been significant analysis of revenue adequacy in markets for electricity transmission capacity \cite{harvey1996transmission}.  Interestingly, it was proven that when accounting for alternating current (AC) power flows in electricity markets, which lead to nonlinear, non-convex continuous optimal power flow problems, revenue adequacy cannot be guaranteed even for simple networks \cite{lesieutre2005convexity}.  In the context of gas pipelines, revenue adequacy was shown to be guaranteed for optimal flow allocation problems provided that zero flow solutions are feasible \cite{rudkevich2017hicss}.  While strict mathematical analysis of revenue adequacy is out of the scope of this study, we discuss several implications that can be drawn from the optimality conditions above.

First, observe that the Lagrange multipliers \(\lambda_j^{NG}\) and \(\lambda_j^{H_2}\) corresponding to constraints \eqref{eq:ngflowbalance} and \eqref{eq:h2flowbalance} are the incremental increases in the objective function value that arise due to incremental increases in withdrawal of natural gas or hydrogen at physical node $j$, respectively.  This means that, in the case that regularity conditions are satisfied, \(\lambda_j^{NG}\) and \(\lambda_j^{H_2}\) can be interpreted as nodal marginal values or prices of natural gas and hydrogen, respectively, with units listed in Table \ref{tab:multipliers_defs} as \$/kg NG and \$/kg $H_2$.  Based on this definition, and to be consistent with the continuity relation \eqref{eq:continuity}, the values \(\lambda_j^{NG}\) and \(\lambda_j^{H_2}\) are associated with the gas at a node after any mixing.  We can compute the value $\lambda_j$ [\$/kg] of the mixture using a linear combination of these prices based on the mass fraction, as well as the locational value of energy $\lambda_j^e$ [\$/MJ] of the blend as
\begin{subequations}
\begin{align} 
& \lambda_j=(1-\gamma_j)\lambda_j^{NG}+\gamma_j\lambda_j^{H_2}, \label{eq:mixture_price} \\
& \lambda_j^e=\frac{(1-\gamma_j)\lambda_j^{NG}+\gamma_j\lambda_j^{H_2}}{ R_{H_{2}}  \gamma_{j(m)} + R_{NG} (1\!-\gamma_{j(m)})}=\frac{\lambda_j}{R(\gamma_j)}, \label{eq:energy_price}
\end{align}
\end{subequations}
where the locational calorific value of the blend $R(\gamma_j)$ is as defined in equation \eqref{eq:energy_conversion}.  The locational values $\lambda_j^{NG}$, $\lambda_j^{H_2}$, $\lambda_j$, and $\lambda_j^e$ correspond to physical nodes $j\in\mathcal{V}$, rather than gNodes $m\in\mathcal{G}$, because the flow balance constraints \eqref{eq:flowbalance} are enforced for physical nodes and we assume no-cost exchange of commodities between market participants co-located at the same physical node.  

Inspecting the first derivative condition \eqref{eq:ldc3} taken with respect to \(d_m\), we can re-arrange this condition and simplify using equation \eqref{eq:energy_price} to obtain an expression for the locational price charged by the market administrator for delivery of energy to a gNode $m$:

\begin{equation} \label{eq:derive_decomposition1}
\lambda_{j(m)}^e  \!=\!  c_m^{d} \!+\! \chi_m^l  \!-\! \chi_m^u \!-\! \chi_m^f  + \frac{c_m^{CO_2}\gamma_{j(m)}}{R(\gamma_{j(m)}) }\!\cdot\!\frac{R_{H_2}}{R_{NG}} \zeta.
\end{equation}
The equation \eqref{eq:derive_decomposition1} can be decomposed as
\begin{subequations} \label{eq:derive_decomposition2}
\begin{align}
\lambda_{j(m)}^e & = \lambda_m^c + \lambda_m^d, \label{eq:energy_price_decomp} \\
\lambda_m^c & =  c_m^{d} + \chi_m^l  - \chi_m^u - \chi_m^f, \label{eq:transportation_component} \\
\lambda_m^d & = \frac{c_m^{CO_2}\gamma_{j(m)}}{R(\gamma_{j(m)}) }\!\cdot\!\frac{R_{H_2}}{R_{NG}} \zeta, \label{eq:decarb_premium}
\end{align}
\end{subequations}
where $\lambda_m^c$ [\$/MJ] is the price component that arises due to the objective component $J_{MR}$ that aims to maximize revenue of the transportation system operator, and $\lambda_m^d$ [\$/MJ] is the \emph{decarbonization premium}.  Observe that the equations  \eqref{eq:transportation_component} for the consumer price component for maximizing revenue $\lambda_m^c$  can be compared with equation (27b) in \cite{rudkevich2017hicss}, which examined locational marginal pricing for homogeneous natural gas pipelines.  The equations are similar, except for an additional term $\chi_m^f$ related to the baseline flow $\bar{g}_m$, which was not considered in the earlier study.

In the case that \(c_{m}^{CO_2} = \$0\)/kg CO\(_2\) and when the gas consumer at \(m\in\mathcal{G}_d\) is a marginal consumer, meaning that neither constraint in equation \eqref{eq:demand_opt} is binding, then the corresponding complementary slackness conditions \eqref{eq:csc11} and \eqref{eq:csc12} imply that  \(\chi_m^l = \chi_m^u = 0\), and equation \eqref{eq:derive_decomposition2} is reduced to \(\lambda_{j(m)}^e = c_m^{d}\).  In the case of congestion arising due to buyer nominations exceeding pipeline capacity, and in the absence of carbon emissions mitigation incentives, the nodal cleared price of energy [\$/MJ] at gNode $m\in\mathcal{G}_d$ is equivalent to the bid price for energy provided by the marginal consumer at that location. When a pipe is at physical pipe capacity, one of the bounds in the set of constraints \eqref{eq:nodallimits} or \eqref{eq:complimits} on pressure, mass fraction, or compression is binding at each end of the pipe, and the cleared market price would reflect such congestion. In the case that $\gamma_j\equiv 0$ system-wide, we have $\lambda_m^d\equiv 0$ for all $m\in\mathcal{G}_d$ in \eqref{eq:decarb_premium}.  In this sense, the price formation done by solving problem \eqref{eq:optimization_formulation} is consistent, because the locational premium on energy prices due to the avoided carbon emissions as given in equation \eqref{eq:carbon_offset} and added with price $c_m^{CO_2}$ in the objective \eqref{eq:obj} has physically meaningful dependence on hydrogen mass fraction, calorific values, molecular weights, and pipeline network capacity.

Now, observe that because all the parameters in equation \eqref{eq:second_order_flow}, namely $f_{ij}$, $L_{ij}$, $D_{ij}$, $A_{ij}$, and the concentration-dependent wave speed $(\gamma_{ij}a_{H_2}^2 + (1-\gamma_{ij})a_{NG}^2)$, are positive, it follows that $\mu_{ij}\cdot\mathrm{sign}(\phi_{ij})<0$.  As a result, $\mathrm{sign}(\mu_{ij}) = -\mathrm{sign}(\phi_{ij})$. Then equation \eqref{eq:ldc4} can be applied to derive
\begin{align}
& \mathrm{sign}\bigl[(\lambda_j^{NG}(1-\gamma_{ij})  + \lambda_j^{H_2}\gamma_{ij}) - (\lambda_i^{NG}(1-\gamma_i) + \lambda_i^{H_2}\gamma_{i})\bigr] \nonumber \\ & \qquad \qquad   = \mathrm{sign}\bigl[(\lambda_j^{NG}(1-\gamma_{ij}) - \lambda_i^{NG}(1-\gamma_i)) + (\lambda_j^{H_2}\gamma_{ij} - \lambda_i^{H_2}\gamma_{i})\bigr] \nonumber \\ 
& \qquad \qquad  = \mathrm{sign}\Bigl[-2\mu_{ij} \frac{f_{ij}L_{ij}}{D_{ij}A_{ij}^2} (\gamma_{ij}a_{H_2}^2 + (1-\gamma_{ij})a_{NG}^2)\left|\phi_{ij}\right|\Bigr] \nonumber \\
& \qquad \qquad = \mathrm{sign}(\phi_{ij}). \label{eq:cash_flow}
\end{align}
In the situation when $\partial_j^+=\{i\}$ and $s_{j(m)}^{NG}= s_{j(m)}^{H_2}\equiv0$ for any $m\in\partial_j^g$, i.e., the only flow incoming to node $j$ is along pipe $(i,j)$, we would have $\gamma_{ij}=\gamma_j$, which would lead to
\begin{equation} \label{eq:cash_flow_mod}
\mathrm{sign}(\lambda_j-\lambda_i) = \mathrm{sign}\bigl[(\lambda_j^{NG}(1-\gamma_{j})  + \lambda_j^{H_2}\gamma_{j}) - (\lambda_i^{NG}(1-\gamma_i) + \lambda_i^{H_2}\gamma_{i})\bigr] = \mathrm{sign}(\phi_{ij}),  
\end{equation}
where we use equation \eqref{eq:mixture_price} for the first equality and equation \eqref{eq:cash_flow} for the second equality.  Equation \eqref{eq:cash_flow_mod} states that the nodal price of gas at the $j$ end of a pipe $(i,j)\in\mathcal{E}$ is higher than at the sending end $i$ if the physical flow is from $i$ to $j$. This result is the intuitive and empirically observable notion that, in the optimal flow allocation solution, gas flows from locations with low price to locations with high price in the pipeline system.  
However, equation \eqref{eq:cash_flow_mod} cannot be guaranteed to hold in the general setting of inhomogeous blending of hydrogen in a pipeline network unless it has a tree structure with any nodal gas injections having lower hydrogen mass fraction than flow from the incoming pipe.  It follows that the marginal price of gas does not monotonically increase in the direction of gas flow in the optimal solution, and we verify this non-monotonicity of dual variables in a case study in Section \ref{sec:baseline}.  Such non-monotonicity is consistent with the non-monotonicity observed in physical variables in previous studies on hydrogen blending \cite{baker2023transitions,zlotnik2023review}, and arises here due to non-monotone changes in hydrogen concentration along pipeline network paths.

Finally, we can observe that the optimal solution to problem \eqref{eq:optimization_formulation} can be used by the pipeline market administrator to collect and appropriately distribute decarbonization incentives in a manner that is consistent with customer participation, pipeline network structure, and the physics of energy flow.  Specifically, suppose that the market administrator collects all of the incentives as defined in the carbon emissions mitigation objective function component $J_{CEM}$ in equation \eqref{eq:obj_components_CEM}.  Then if $\lambda_m^d$ is the decarbonization premium as in equation \eqref{eq:decarb_premium}, and $d_m$ is the physical flow at hydrogen mass fraction $\gamma_{j(m)}$ to the consumer at gNode $m\in\mathcal{G}_d$, the equitable pass-through distribution of the incentives is accomplished by distributing a credit $\tau_m$ [\$/s] defined by 
\begin{equation} \label{eq:pass_through_credit}
\tau_m=\lambda_m^d\cdot d_m \cdot \left( R_{H_{2}}  \gamma_{j(m)} + R_{NG} (1\!-\gamma_{j(m)}) \right)
\end{equation}
to the consumer at gNode $m\in\mathcal{G}_d$.  We define the total amount of distributed pass-through credits $D_{PTC}$ [\$/s] as 
\begin{equation} \label{eq:pass_through_total}
D_{PTC}=\sum_{m\in\mathcal{G}_d} \tau_m.
\end{equation}
We suppose that, if the problem \eqref{eq:optimization_formulation} is solved using a primal-dual interior point method such that there is no dual feasibility or optimality gap, the multipliers $\lambda_j^{NG}$ and $\lambda_j^{H_2}$ correctly quantify sensitivities to incremental nodal consumptions of these commodities.  In that case, $\lambda_m^d$ would also be consistently quantified for all $m\in\mathcal{G}$ such that $D_{PTC}\equiv J_{CEM}$, that is, the carbon emissions mitigation incentives collected by the administrator are sufficient to distribute all pass-through credits according to the nodal decarbonization premium.   We will examine this property in the case studies in the following section. 

\section{Computational Case Studies}  \label{sec:case_study}

We demonstrate the use of the proposed optimization formulation for operating a pipeline market mechanism for natural gas and hydrogen transport using an 8-node network and a 40-node network. The problem \eqref{eq:optimization_formulation} is implemented in the Julia programming language v1.7.3 using the JuMP optimization toolkit v1.7.0 \cite{dunning2017jump}. Our implementation makes use of the general purpose interior point solver IPOPT v1.6.2 \cite{wachter2006implementation}, a large-scale nonlinear optimization software. The case studies are evaluated using an Apple M1 chip with 8 GB of RAM, with solve time of about 60 milliseconds for both networks.  The code and data sets developed for the computations in this study is available in open source \cite{phydras2024github}.

\subsection{Re-scaling, Non-Dimensionalization, and Parameters} \label{sec:rescale}

Prior to solving problem \eqref{eq:optimization_formulation}, we non-dimensionalize the physical variables in the governing equations in order to avoid numerical issues \cite{srinivasan2022numerical}. We re-scale equation \eqref{eq:pipe_equations} because the squared wave speed \(V\) in blended gas is much larger than other parameters in the equation. Given the transformations \(\bar{P} = P / P_0\), \(\bar{L} = L/ l_0\), \(\bar{D} = D/ l_0\), \(\bar{A} = A/ A_0\), and \(\bar{\phi} = \phi/ \phi_0 = \phi / (\rho_0 u_0 A_0)\), equation \eqref{eq:pipe_equations} becomes 
\begin{subequations} \label{eq:pipe_equations_dimensionless}
\begin{align}
    \bar{P}_{i}^{2} - \bar{P}_{j}^{2} = \frac{f_{ij}\bar{L}_{ij}}{\bar{D}_{ij}\bar{A}_{ij}^{2}} \bar{V}_{ij}  \bar{\phi_{ij}} \left|\bar{\phi}_{ij} \right| \cdot \frac{u_0^2}{a_0^2},  \quad \forall (i,j) \in \mathcal{E}, \label{eq:weymouth_dimensionless} \\ 
    \bar{V}_{ij} \triangleq \frac{V_{ij}(\gamma_{ij})}{a_0^2},  \quad \forall (i,j) \in \mathcal{E}. \label{eq:soundspeed_dimensionless}
\end{align}
\end{subequations}
The nominal length, area, density, and velocity used for both network models are \(l_0 = 5000\) m, \(A_0 = 1\) m$^2$, \(\rho_0 = P_0 / a_0^2\), and \(u_0 = \lceil a_0 \rceil / 300\), where \(a_0\) is the geometric mean of wave speeds in the gases.  We compute \(a_0\) as \(a_0 = \sqrt{a_{NG} \cdot a_{H_2}} \)  where $a_{NG}$ and $a_{H_2}$ are the wave speeds in NG and H$_2$, obtained by \(a_{NG} = \sqrt{RT/M_{NG}}\) and \(a_{H_2} = \sqrt{RT/M_{H_2}}\), respectively. Here, \(R = 8.314\) J/mol/K is the universal gas constant, and \(M_{NG} = 0.01737\) kg/mol and \(M_{H_2} = 0.002016\) kg/mol are molecular masses of NG and H$_2$. Using these parameters results in a nominal wave speed of \(a_0 = 635.06\) m/s used for re-scaling. The parameter values used in our study are summarized in Table \ref{tab:parameter_allnode}.

\begin{table}[!h]
\centering \normalsize
\begin{tabular}{|p{3cm}|p{5cm}||p{3cm}|p{3cm}|}
  \hline
  Parameters  & Values & Parameters  & Values \\ 
  \hline
  \hline
  $G_{H_2}, G_{NG}, G_{nom}$  & 0.0696, 0.6, 0.574 & $P_{j}^{min}, P_{j}^{max}$   & 3 MPa, 6 MPa \\
  $\kappa_{H_2}, \kappa_{NG}, \kappa_{nom}$  & 1.405, 1.303, 1.308 & $\gamma_{j}^{min}, \gamma_{j}^{max}$   & 0, 0.1  \\
  $a_{H_2}, a_{NG}$  & 1090 m/s, 370 m/s & $s_{m}^{H_2,max}, s_{m}^{NG,max}$   & inf kg/s \\
  $R_{H_2}, R_{NG}$  & 141.8 MJ/kg, 44.2 MJ/kg & $c_{m}^{d}$   & \$0.019/MJ  \\
  $T$   & 288.7 K & $c_{m}^{H_2}, c_{m}^{NG}$   & \$0.8/kg, \$0.2/kg \\
  $\zeta$   & 44/16 & & \\
  \hline
\end{tabular}
\vspace{-1.5ex}
\caption{Physical parameter values for 8-node \& 40-node networks.} \label{tab:parameter_allnode}
\end{table}

The natural gas and hydrogen prices listed in Table \ref{tab:parameter_allnode} are intended to represent a baseline for prevailing market conditions, in which hydrogen is more expensive per unit mass as well as per unit energy, and these factors are 4 and $\approx\!\!1.25$, respectively, as can be computed using the data.  In reality, the prevailing price of energy in the form of hydrogen is typically significantly greater than these assumptions.  Natural gas exhibits volatile price fluctuations \cite{macavoy2008natural}, and a comprehensive analysis of prevailing prices, techno-economic assessments, and forecasts for hydrogen trade values can be found in a recent study \cite{franzmann2023green}.  Our aim here is to develop a methodology for optimal scheduling and pricing of natural gas and hydrogen transport over a pipeline network given any specified market conditions, and further review of prevailing or future prices is beyond our intended scope.

\begin{figure}[h!]
\centering
\includegraphics[width=.8\columnwidth]{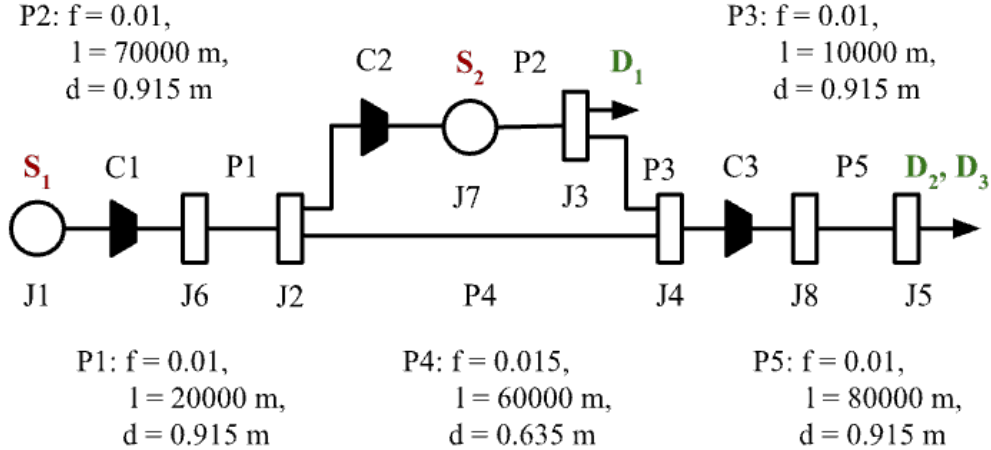} 
\caption{8-node test network schematic, with physical characteristics. The network consists of three compressors, one natural gas supplier and one hydrogen supplier, and three flexible offtakers. \label{fig:8node}}
\end{figure}

\subsection{8-Node Network Case Study} \label{sec:case1}

Our first case study applies problem \eqref{eq:optimization_formulation} to the 8-node network shown in Figure \ref{fig:8node}.  This case study builds on our previous  sensitivity analysis for this same system \cite{sodwatana2023h2blend}.  The purpose of examining this small test network is because solutions to problem \eqref{eq:optimization_formulation} can undergo qualitative transitions of interest while being simple enough to interpret.  We suppose that the network has two suppliers that provide natural gas and hydrogen, respectively, without supply quantity constraints.  The natural gas supplier S1 is located at physical node J1, and the hydrogen supplier S2 is located downstream at node J7. There is an offtaker D1 at node J3 and two offtaker gNodes, D2 and D3, located at physical node J5. The withdrawals of the three offtakers are flexible, following equation \eqref{eq:demand_opt}, so that their actual energy consumption is constrained at \(g_m^{max} = 2000\) MJ/s and is determined by the market solution produced by solving problem \eqref{eq:optimization_formulation}. In addition, there are three compressors, each with a maximum pressure boost ratio of 1.4. The maximum allowable hydrogen injection is 10\% by mass. The network has one slack node (J1) with nominal pressure \(P_0 = 4\) MPa. The complete operational constraints and pipeline characteristics are summarized in Figure \ref{fig:8node}, and parameters are given in Table \ref{tab:parameter_allnode}. We consider two market scenarios for the 8-node network.

\begin{table}[!h]
\centering
\begin{tabular}{|p{3cm}|c|c|c|c|c|}
 \hline
 gNode & S1 (NG) & S2 (H\(_2\)) & D1 & D2 & D3 \\
 Physical node & J1 & J7 & J3 & J5 & J5 \\
 \hline
 \hline
 NG flow [kg/s] & 135 & 0 & 45 & 45 & 45\\
 H\(_2\) flow [kg/s] & 0 & 0 & 0 & 0 & 0\\
 Total flow [kg/s] & 135 & 0 & 45 & 45 & 45\\
 Energy [MJ/s] & 6000 & 0 & 2000 & 2000 & 2000\\
 Pressure [MPa] & 4.00 & 3.84 & 3.50 & 3.14 & 3.14\\
 CI [kgCO$_2$/MJ] & - & - & 0.0622 & 0.0622 & 0.0622\\
 \(\gamma_j\) [-] & 0 & 0 & 0 & 0 & 0 \\
 \(\lambda_{j}^{NG}\) [\$/kg] & 0.20  & 0.20 & 0.20 & 0.20 & 0.20\\
 \(\lambda_{j}^{H_2}\) [\$/kg] & -  & - & -  & - & - \\
 \(\lambda_{j}^e\) [\$/MJ] & 0.0045 & 0.0045 & 0.0045 & 0.0045 & 0.0045\\
 \(\lambda_m^d\) [\$/MJ] & - & - & 0 & 0 & 0\\
 \hline
\end{tabular}
\caption{Optimal market solution for 8-node network case study Scenario 1, with $c_m^{CO_2}=0$ $\forall m\in\mathcal{G}$.  The objective value is \(J_{EV}=86.85\), with components \(J_{MR}=86.85\), \(J_{CEM}=0\), and \(J_{GC}=0\), in units of [\$/s]. \label{tab:case1a}}
\end{table}

In Scenario 1, we suppose that the system-wide carbon emissions mitigation incentive is \(c_{m}^{CO_2} = \$0\)/kg CO\(_2\). The optimal solution to problem \eqref{eq:optimization_formulation} in this scenario is appropriate in the sense that the physical and market solutions are consistent with the pipeline network capacity and market conditions defined by the price-quantity bids of market participants. Specifically, the pipeline operator is able to schedule flows such that all three customers receive amount of energy specified in their quantity bids $g_m^{\max}$ in the form of natural gas. The cleared market price for the delivered gas is \$0.0045/MJ, which is equivalent to the supply price of natural gas, indicating no pipeline congestion. Because there were no CO$_2$ mitigation incentives, it was not profitable for hydrogen to be injected into the pipeline network, and the decarbonization premium \(\lambda_m^d\) is zero at all withdrawal nodes. The objective function value in Scenario 1 is \(J_{EV}^1 = J_{MR} + J_{GC} = \$86.85/s\) and the total carbon dioxide emitted from natural gas consumption is 373 kg/s.  

\begin{table}[!h]
\centering
\begin{tabular}{|p{3cm}|c|c|c|c|c|}
 \hline
 gNode & S1 (NG) & S2 (H\(_2\)) & D1 & D2 & D3 \\
 Physical node & J1 & J7 & J3 & J5 & J5 \\
 \hline
 \hline
 NG flow [kg/s] & 110 & 0 & 33.3 & 38.5 & 38.5\\
 H\(_2\) flow [kg/s] & 0 & 7.7 & 3.7 & 2 & 2\\
 Total flow [kg/s] & 110 & 7.7 & 37 & 40.5 & 40.5\\
 Energy [MJ/s] & 4908 & 1092 & 2000 & 2000 & 2000\\
 Pressure [MPa] & 4.00 & 3.89 & 3.52 & 3.11 & 3.11\\
 CI [kgCO$_2$/MJ] & - & - & 0.0458 & 0.0527 & 0.0527\\
 \(\gamma_j\) [-] & 0 & 0.10 & 0.10 & 0.05 & 0.05 \\
 \(\lambda_{j}^{NG}\) [\$/kg] & 0.20  & 0.20 & 0.20 & 0.18 & 0.18\\
 \(\lambda_{j}^{H_2}\) [\$/kg] & -  & 0.80 & 0.94  & 1.07 & 1.07 \\
 \(\lambda_{j}^e\) [\$/MJ] & 0.0045 & 0.0048 & 0.0050 & 0.0046 & 0.0046\\
 \(\lambda_m^d\) [\$/MJ] & - & - & 9.0e-4 & 5.2e-4 & 5.2e-4 \\
 \hline
\end{tabular}
\caption{Optimal market solution for 8-node network case study Scenario 2, with $c_m^{CO_2}=\$0.055$ $\forall m\in\mathcal{G}$. The objective value is \(J_{EV}=89.46\), with components \(J_{MR}=85.59\), \(J_{CEM}=3.87\), and \(J_{GC}=0\), in [\$/s]. \label{tab:case1b}}
\end{table}

In Scenario 2, we introduce an emissions offset incentive of \(c_{m}^{CO_2} = \$0.055\)/kg CO\(_2\) for all $m\in\mathcal{G}$, which is approximately \$50 per U.S. ton. All three consumers still receive their requested amount of energy. However, the total amount of natural gas injected into the system decreases by 25 kg/s, from 135 kg/s to 110 kg/s, which is a 23\% decrease, and the equivalent energy is replaced with 7.7 kg/s of hydrogen. Consumer D1, which is immediately downstream of the hydrogen injection point, receives a 10\% hydrogen blend, while consumers D2 and D3 receive a 5\% hydrogen blend. Examining the market solutions, we see that the locational market price for the energy delivered to gNode D1 increases to \$0.0050/MJ and to \$0.0046/MJ at D2 and D3, which are 11\% and 2\% increases compared to Scenario 1, respectively. The differences in price between the two physical nodes reflect the difference in the locational values of energy and of the decarbonization premium paid by the end-users. The decarbonization premium \(\lambda^d\) at D1 is \$0.00090/MJ while the decarbonization premium at D2 and D3 is \$0.00052/MJ, and these values are 18\% and 11.3\% of the locational prices, respectively. Overall, the economic value produced by the pipeline operator in Scenario 2 increases to \(J_{EV}^2 = \$89.46/s\), which is about 3\% higher than $J_{EV}^1$ in Scenario 1,  while the total carbon dioxide emitted is 303 kg/s, which is a 23\% reduction with respect to Scenario 1.  The carbon intensity of energy delivered at D1 and at D2 and D3 is reduced by 35\% and 18\%, respectively. The solutions obtained by solving \eqref{eq:optimization_formulation} are shown in Table \ref{tab:case1a} for Scenario 1 and Table \ref{tab:case1b} for Scenario 2. Note that \(\lambda_j(m)\), the marginal price of the blended gas, is presented in units of [\$/MJ] by dividing through by \(R(\gamma_{j(m)})\). The binding constraints in Scenario 2 that limit the amount of H\(_2\) blended are the blend concentration constraints of 10\% by mass between the H\(_2\) injection site (S2) and the merging node (J4). If this blend constraint bound were increased, the system could increase the total objective function value by increasing hydrogen blending.

\subsection{40-Node Network Case Studies} \label{sec:case2}

We now examine solutions to problem \eqref{eq:optimization_formulation} when applied to a more complex pipeline network model in order to demonstrate that the implementation can be generalized to arbitrary network topologies and a variety of market conditions.  We make use of a 40-node pipeline network model for this purpose, and we represent the values of physical and market solution variables using color in illustrations of the network, which are more interpretable than large tables.  We first examine the solution for a baseline scenario that includes some hydrogen blending, and then modify it in specific ways to examine how well-intended policies may result in unintended  outcomes under certain market conditions.

\begin{figure*}[!ht]
\includegraphics[width=\textwidth]{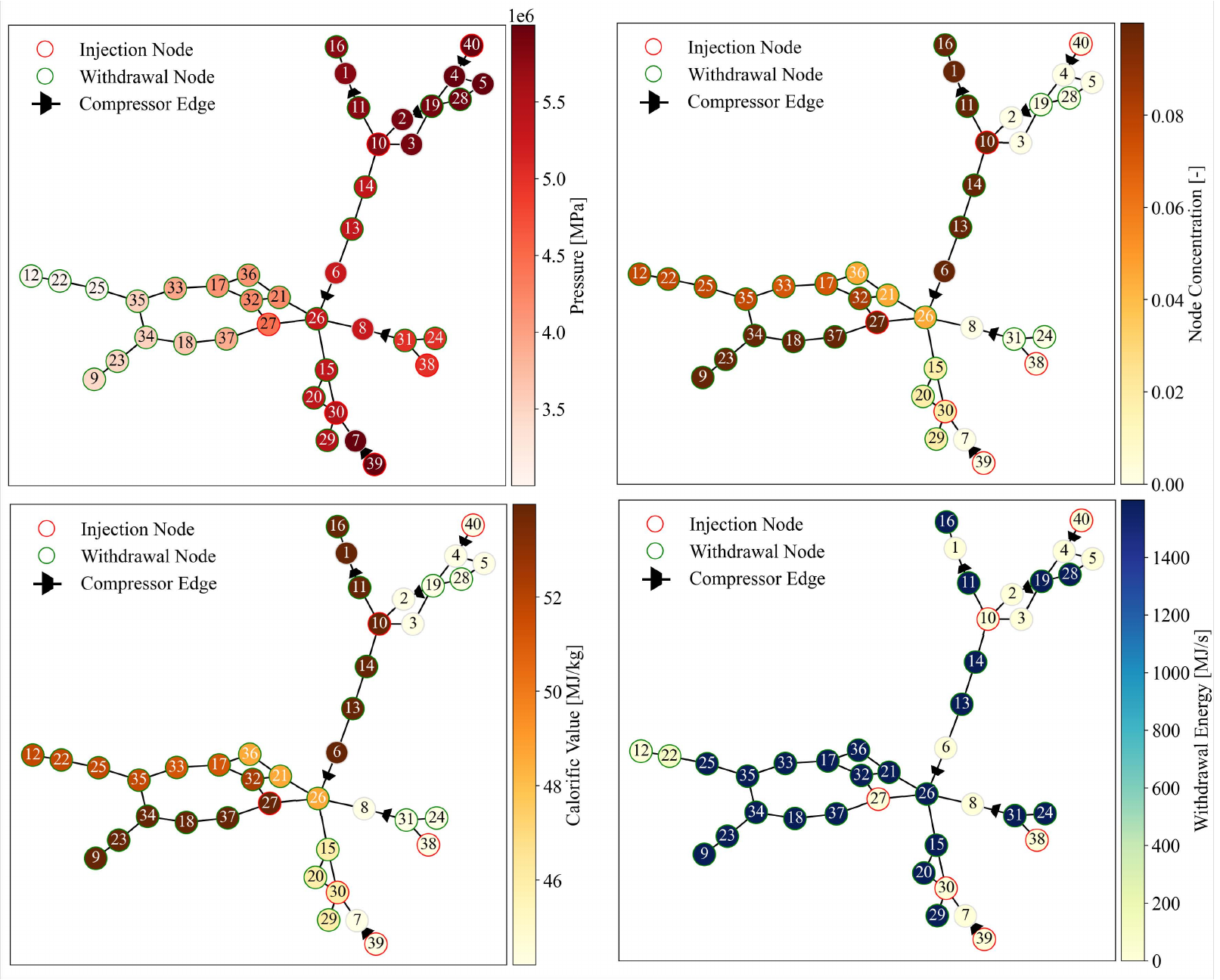} 
\centering
\caption{Physical solution for the example Scenario using the 40-node network. From left to right: top row shows pressure and nodal hydrogen concentration; bottom row shows calorific value and withdrawal energy. \label{fig:example_output_physical}}
\end{figure*}

The 40-node network case study is a modification of the GasLib-40 network \cite{schmidt2017gaslib}, with one slack node (node 38) with a nominal pressure \(P_0 = 5\) MPa, and 26 physical withdrawal nodes. Energy withdrawal by all customers is flexible, following equation \eqref{eq:demand_opt}. The bid prices and quantities of all consumers are set to \(c_m^{d}\) = \$0.019/MJ and \(g_m^{max}\) = 1600 MJ/s. We consider a system-wide carbon emissions mitigation incentive price of \(c_m^{CO_2}\) = \$0.055/kg for all $m\in\mathcal{G}$. We use the model parameters listed in Table \ref{tab:parameter_allnode}, as for the 8-node network case study.  The network topology, including locations of injection/supply nodes, withdrawal nodes, and compressors is indicated in the plots of results in Figures \ref{fig:example_output_physical} and \ref{fig:example_output_market}.

\subsubsection{Baseline Scenario} \label{sec:baseline}

\begin{figure*}[!t]
\includegraphics[width=\textwidth]{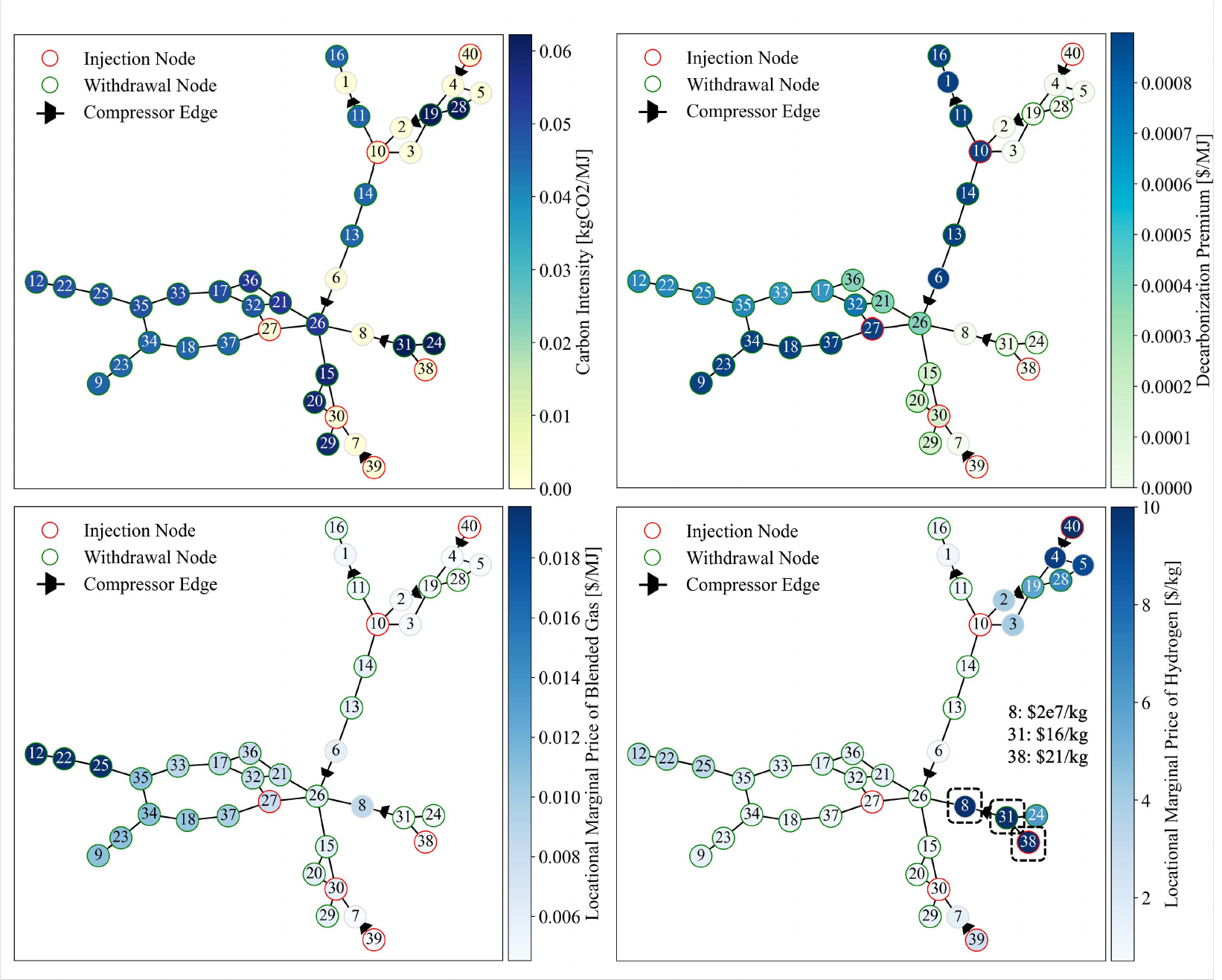} 
\centering
\caption{Market solution for the example Scenario using the 40-node network. From left to right: top row - carbon intensity and decarbonization premium; bottom row - LMP of the blended gas and LMP of hydrogen. For graphical clarity, we limit the values of the nodes in dashed boxes in the bottom right plot to be 10, with actual values listed above. \label{fig:example_output_market}}
\vspace{-1ex}
\end{figure*}

In this baseline scenario, we suppose there are three physical nodes (38, 39 and 40) at which suppliers inject natural gas with unlimited supply, and three separate physical nodes (10, 27 and 30) with suppliers that inject hydrogen.  We show visualizations of the optimal physical and market solutions to problem \eqref{eq:optimization_formulation} throughout the network in Figures \ref{fig:example_output_physical} and \ref{fig:example_output_market}, respectively.  Inspecting these results, we observe that the optimal solution to problem \eqref{eq:optimization_formulation} behaves as expected in both the physical (primal) and market (dual) variables.

To examine the physical solution, consider first the nodal pressures as illustrated in the top left panel in Figure \ref{fig:example_output_physical}. The pressure is highest at the natural gas supply nodes 38, 39, and 40 in the lower and upper right parts of the map, and is lowest at nodes 9 and 12, which are at the far downstream ends of the pipeline at the left part of the figure.  Next, observe the nodal hydrogen concentration as illustrated in the top right panel in Figure \ref{fig:example_output_physical}.  We see that the hydrogen mass fraction is near zero at the natural gas supply nodes 38, 39, and 40 and downstream from them until the flow reaches nodes 10 or 27, which have substantial hydrogen injections up to 10\% mass fraction.  Because of the looped topology of the network, there are multiple nodes, such as nodes 26, 32, and 35, which show intermediate mass fractions that do not bind at the allowable maximum $\gamma_j^{\max}$.  Inspecting the calorific value as illustrated in the bottom left panel in Figure \ref{fig:example_output_physical}, it is clear that there is a one-to-one correspondence to the nodal hydrogen concentration results in the top right panel.  Note that this calorific value is shown as a per-mass value rather than a per-volume value.  We expect lower volumetric calorific value with increasing hydrogen injection, corresponding to higher flow velocities to transport the same energy content. Finally, we observe in the bottom right panel in Figure \ref{fig:example_output_physical} that 24 of the 26 physical withdrawal nodes that provide a price-quantity bid into the market mechanism receive the 1600 MJ/s allocation that is the upper bound value for optimized energy withdrawal.  Nodes 12 and 22, however, are not allocated their quantity bids, likely because they are furthest downstream from the supply locations.  These are the marginal consumers, and so the LMPs at nodes 12 and 22 are at their price bids of \$0.019/MJ as seen at bottom left in Figure \ref{fig:example_output_market}. Based on these results, we consider that the solution to  optimization problem \eqref{eq:optimization_formulation} behaves in a physically realistic manner, i.e., that pressure decreases in the direction of flow, and mass balance of the gas constituents is appropriately conserved.

Inspecting the market solution, we see that the nodal carbon intensity as illustrated in the top left panel in Figure \ref{fig:example_output_market} is highest at those locations where the nodal hydrogen concentration is lowest, seen in the top right panel in Figure \ref{fig:example_output_physical}.  We can also see that the decarbonization premium paid by consumers shown in top right panel in Figure \ref{fig:example_output_market} is highest at those locations with the greatest nodal hydrogen concentration, which is shown in the top right panel in Figure \ref{fig:example_output_physical}.  This is in accordance with the decarbonization premium as defined in equation \eqref{eq:decarb_premium}.  Next, examining the bottom left panel in Figure \ref{fig:example_output_market}, we see that the locational price of the gas mixture is lowest at the supply nodes (nodes 38, 39, and 40) and highest at the downstream nodes (nodes 9 and 12), and by comparing with the top left panel in Figure \ref{fig:example_output_physical} we see that price generally increases as pressure decreases in the direction of flow. 
However, this is not strictly the case, as can be seen by inspecting the solution at nodes 6, 8, and 26.  The pressures are consistent with gas flowing from nodes 6 and 8 into 26, and the nodes have hydrogen concentrations of 0.1, 0.0, and 0.0455, respectively, and energy LMPs of \$0.00588/MJ, \$0.00868/MJ, and \$0.00603/MJ, respectively.  The decrease in the price of energy along edge (8,26) is explained by the blending of gas with higher hydrogen concentration coming in from edge (6,26).
Next, inspecting the bottom right panel in Figure \ref{fig:example_output_market}, and comparing with the top left panel in Figure \ref{fig:example_output_physical}, we can see that the locational price of hydrogen can decrease as pressure decreases in the direction of flow.  This can be explained by the need to provide incentives without which hydrogen suppliers would not participate in the market.  In fact, as indicated in the bottom right panel in Figure \ref{fig:example_output_market}, the nodes 38, 31, and 8 that are upstream of any hydrogen injection have the highest hydrogen values.  Therefore we propose that solving the problem \eqref{eq:optimization_formulation} for a given network and market structure indicates the ideal locations (those with highest $\lambda_j^{H_2}$) at which adding incremental hydrogen suppliers will maximize mitigation of carbon emissions.



\begin{figure*}[!h]
\includegraphics[width=\textwidth]{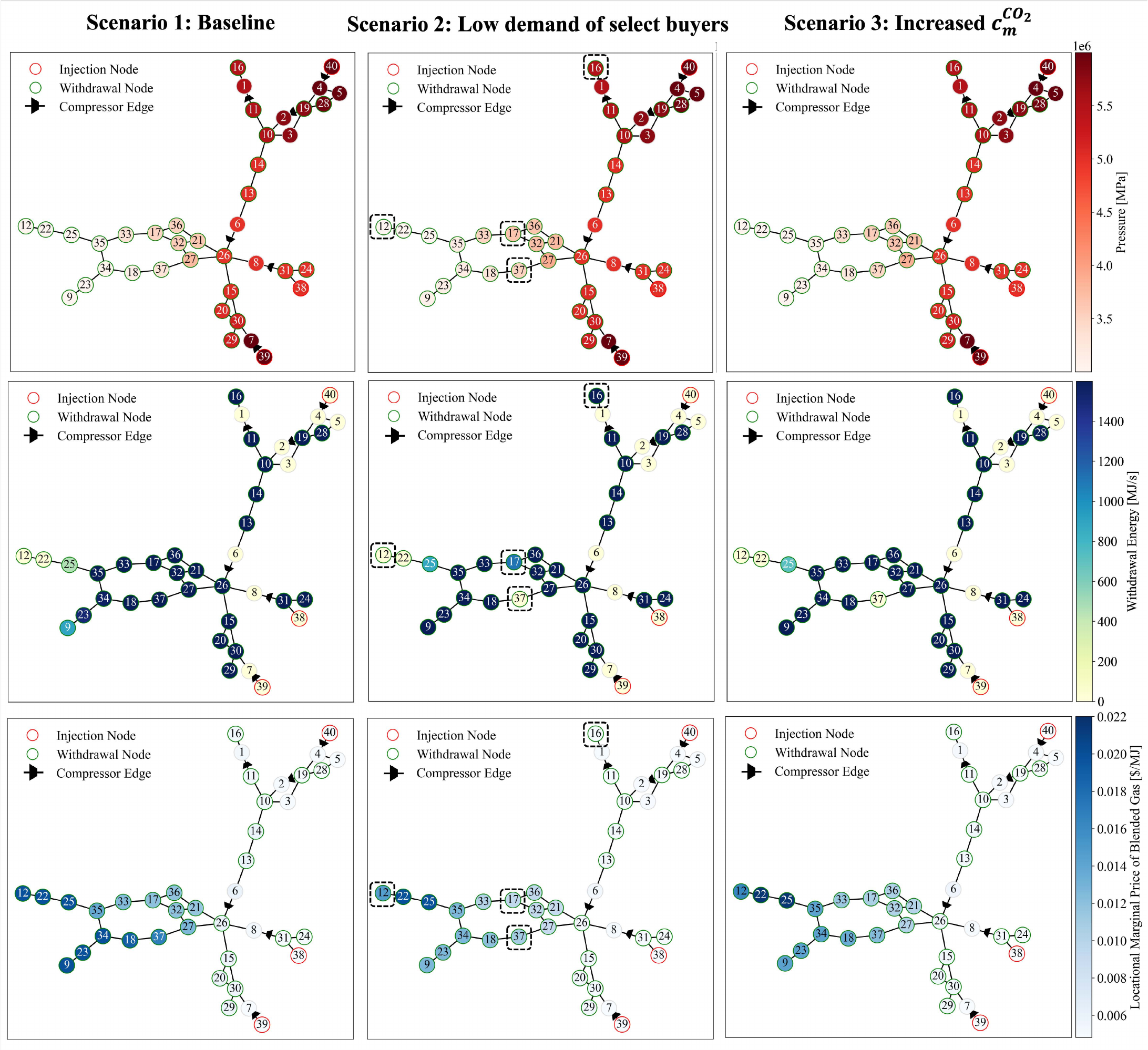} 
\centering
\caption{Results from the three scenarios of the counter-intuitive case with the 40-node test network. Each column corresponds to a scenario: the left column corresponding to Scenario 1; the middle column refers to Scenario 2, where \(c_m^d\) of four gNodes represented in the dashed box decreased to \$0.0085/MJ; the right column is for Scenario 3, where we increase \(c_m^{CO_2}\) system-wide to \$0.155/kg CO\(_2\) for all $m\in\mathcal{G}$. Top to bottom: pressure [MPa]; optimal withdrawal energy [MW]; locational marginal price of the blended gas [\$/MJ]. The numbers within each node corresponds to the physical node number. Injection nodes are circled in red and consists of physical nodes 38, 39, and 40. The 29 withdrawal nodes are circled in green. \label{fig:40node_summary}}
\end{figure*}

\subsubsection{Counter-intuitive Market Outcomes} \label{sec:case_counter}

Next, we develop a price structure for the 40-node network in which decarbonization incentives cause greater carbon emissions in the optimal market solution.  We suppose that the network has three physical injection nodes (nodes 38, 39 and 40), where each has two gNodes that inject hydrogen and natural gas. Hydrogen injection sites in the previous case are now withdrawal nodes, thus we have a total of 29 withdrawal nodes in this case study. All other parameters are the same as in the baseline scenario.

We start with Scenario 1, in which the price and quantity bids for all offtake gNodes are \(c_m^{d}\) = \$0.019/MJ and \(g_m^{max}\) = 1600 MJ/s with a system-wide carbon mitigation incentive \(c_m^{CO_2}\) fixed at \$0.055/kg for all $m\in\mathcal{G}$. In Scenario 2, we suppose that four offtakers at nodes 12, 16, 17, and 37, indicated using green squares in Fig. \ref{fig:40node_summary}, decrease their price bids for energy by half to \(c_m^{d}\) = \$0.0085/MJ to reflect low demand.  This scenario could represent the change in the market position of a gas-fired power plant at a time with low prevailing prices in the wholesale electricity market. In Scenario 3, we suppose that \(c_m^{CO_2}\) is increased to \$0.155/kg for all $m\in\mathcal{G}$, which represents a scenario with low demand and high decarbonization incentives. We compare the optimal solution to problem \eqref{eq:optimization_formulation} for the three scenarios in Figure \ref{fig:40node_summary}, which shows nodal values for pressure, energy consumption, and the locational marginal price \(\lambda_j\)  of the blended gas.  The aggregate optimal outcomes are compared for the scenarios in addition to the baseline in the previous section in Table \ref{tab:fortynode_results}.  More localized results for the three scenarios are shown in Tables \ref{tab:case2a}, \ref{tab:case2b}, and \ref{tab:case2c}, to compare the solutions at the withdrawal nodes 12, 16, 17, and 37 where the price bids are decreased in Scenarios 2 and 3.

\begin{table}[!h]
\centering
\begin{tabular}{|p{3cm}|c|c|c|c|}
 \hline
 Physical node & 16 & 17 & 37 & 12 \\
 \hline
 \hline
 NG flow [kg/s] & 26.69 & 26.69 & 26.69 & 26.69 \\
 H\(_2\) flow [kg/s] & 2.96 & 2.96 & 2.96 & 2.96 \\
 Total flow [kg/s] & 29.65 & 29.65 & 29.65 & 0.02 \\
 Energy [MJ/s] & 1600 & 1600 & 1600 & 1 \\
 Pressure [MPa] & 5.48 & 3.49 & 3.30 & 3.00 \\
 CI [kgCO$_2$/MJ] & 0.0459 & 0.0459 & 0 & 0 \\
 \(\gamma_j\) [-] & 0.10 & 0.10 & 0.10 & 0.10  \\
 \(\lambda_{j}^{NG}\) [\$/kg] & 0.1965 & 0.5068 & 0.7212 & 7.5579 \\
 \(\lambda_{j}^{H_2}\) [\$/kg] & 1.0507 & 1.8674 & 2.6086 & -57.28 \\
 \(\lambda_{j}^e\) [\$/MJ] & 0.0052 & 0.0119 & 0.0169 & 0.0199 \\
 \(\lambda_m^d\) [\$/MJ] & 8.99e-4 & 8.99e-4 & 8.99e-4 & 8.99e-4 \\
 \hline
\end{tabular}
\caption{Optimal market solution at four withdrawal nodes for Scenario 1 of the 40-node network case study. \label{tab:case2a}}
\vspace{-3ex}
\end{table}

In Scenario 1, the pipeline is constrained by the lower and upper pressure limits, which determine the maximum flow capacity on edges. Hydrogen is blended at the supply nodes at its maximum limit of 10\% mass fraction. In Figure \ref{fig:40node_summary} left top panel, we see that nodes far from injection sites and compressor have significantly lower pressure.  Inspecting the left second row panel, we see that delivery of energy at the downstream nodes 9 and 25 is curtailed below the quantity bids of the consumers at those locations. Comparing the locational price solution in the left bottom panel with the left top panel in Figure \ref{fig:40node_summary}, we see that the cleared market price for the delivered gas increases as pressure decreases in the direction of flow. The objective function value is \(J_{EV} = \$622/s\) and the total carbon dioxide emitted is 1891 kg/s.  The results of the solution in Scenario 1 are summarized for the selected nodes in Table \ref{tab:case2a}.

When the price bids of four offtakers at nodes 12, 16, 17, and 37 are reduced in Scenario 2, the optimal solution to problem \eqref{eq:optimization_formulation} results in hydrogen being blended at its maximum limit of 10\%. Node 16 still receives its quantity bid for energy, at the maximum constraint bound value. Inspecting the pressure solution in the center top panel in Figure \ref{fig:40node_summary}, we see that nodal pressure at node 16 remains high at 5.48 MPa given its proximity to a compressor. However, the pipeline now delivers only 1110 MW to node 17 while node 37 receives no energy at all, which can be justified given the lower prices bid at that node. Node 12 receives no energy delivery as well. The objective function value is \(J_{EV} = \$580/s\) and the total carbon dioxide emitted is 1853 kg/s.  The results of the solution in Scenario 2 are summarized for the selected nodes in Table \ref{tab:case2b}.

\begin{table}[h]
\centering
\begin{tabular}{|p{3cm}|c|c|c|c|}
 \hline
 Physical node & 16 & 17 & 37 & 12 \\
 \hline
 \hline
 NG flow [kg/s] & 26.69 & 18.52 & 0.00 & 0.00 \\
 H\(_2\) flow [kg/s] & 2.96 & 2.06 & 0.00 & 0.00 \\
 Total flow [kg/s] & 29.65 & 20.58 & 0.00 & 0.00 \\
 Energy [MJ/s] & 1600 & 1110 & 0 & 0 \\
 Pressure [MPa] & 5.48 & 3.65 & 3.52 & 3.00 \\
 CI [kgCO$_2$/MJ] & 0.0459 & 0.0459 & 0 & 0 \\
 \(\gamma_j\) [-] & 0.10 & 0.10 & 0.10 & 0.10  \\
 \(\lambda_{j}^{NG}\) [\$/kg] & 0.1873 & 0.4024 & 0.4622 & 3.1e4 \\
 \(\lambda_{j}^{H_2}\) [\$/kg] & 1.0152 & 1.4502 & 1.7292 & -2.8e5 \\
 \(\lambda_{j}^e\) [\$/MJ] & 0.0050 & 0.0094 & 0.0109 & 0.0146 \\
 \(\lambda_m^d\) [\$/MJ] & 8.99e-4 & 8.99e-4 & 8.99e-4 & 8.99e-4 \\
 \hline
\end{tabular}
\caption{Optimal market solution at four withdrawal nodes for Scenario 2 of the 40-node network case study. \label{tab:case2b}}
\vspace{-3ex}
\end{table}

In Scenario 3, we increase the system-wide carbon reduction incentive \(c_m^{CO_2}\) to \$0.155/kg for all $m\in\mathcal{G}$. In the optimal solution to problem \eqref{eq:optimization_formulation}, a pressure of 5.48 MPa is maintained at node J17, at which the consumer withdraws its requested quantity bid of 1600 MW of energy. This quantity is a 45\% increase from the 1110 MW of energy delivered to node 17 in Scenario 2. Meanwhile, nodes 12 and 37 receive no energy. The total energy delivered by the system increases although the hydrogen fraction remains constant at the upper limit of 10\%. Consequently, this change leads to more natural gas consumption and higher carbon dioxide emissions. The objective function value is \(J_{EV} = \$647/s\) and the total carbon dioxide emitted is 1868 kg/s, which is a slight increase with respect to Scenario 2.  The results of the solution in Scenario 3 are summarized for the selected nodes in Table \ref{tab:case2c}.

\begin{table}[h]
\centering
\begin{tabular}{|p{3cm}|c|c|c|c|}
 \hline
 Physical node & 16 & 17 & 37 & 12 \\
 \hline
 \hline
 NG flow [kg/s] & 26.69 & 26.69 & 0.00 & 0.00 \\
 H\(_2\) flow [kg/s] & 2.96 & 2.96 & 0.00 & 0.00 \\
 Total flow [kg/s] & 29.65 & 29.65 & 0.00 & 0.00 \\
 Energy [MJ/s] & 1600 & 1600 & 0 & 0 \\
 Pressure [MPa] & 5.48 & 3.56 & 3.47 & 3.00 \\
 CI [kgCO$_2$/MJ] & 0.0459 & 0.0459 & 0 & 0 \\
 \(\gamma_j\) [-] & 0.10 & 0.10 & 0.10 & 0.10  \\
 \(\lambda_{j(m)}^{NG}\) [\$/kg] & 0.1357 & 0.4486 & 0.5085 & 1.4e5 \\
 \(\lambda_{j(m)}^{H_2}\) [\$/kg] & 1.5110 & 1.6454 & 2.0453 & -1.2e6 \\
 \(\lambda_{j(m)}^e\) [\$/MJ] & 0.0051 & 0.0105 & 0.0123 & 0.0163 \\
 \(\lambda_m^d\) [\$/MJ] & 2.52e-3 & 2.52e-3 & 2.52e-3 & 2.52e-3 \\
 \hline
\end{tabular}
\caption{Optimal market solution at four withdrawal nodes for Scenario 3 of the 40-node network case study. \label{tab:case2c}}
\vspace{-3ex}
\end{table}

Inspecting the values of the decarbonization premium \(\lambda_m^d\) in \eqref{eq:derive_decomposition2}, we observe that increasing the carbon emissions mitigation incentive \(c_m^{CO_2}\) increases the cleared market price for delivered gas. In the analysis in Section \ref{sec:case1}, where energy demand is met entirely by natural gas and the there are no binding constraints on hydrogen injection in the optimal solution, this increase in \(c_m^{CO_2}\) encourages a substitution of natural gas with some hydrogen. However, in this analysis, there are unmet energy demands that result because of unfavorably low bid prices in Scenario 2 rather than any operational constraints. Increasing \(c_m^{CO_2}\) in Scenario 3 results in the optimal market solution allocating energy to customers with low price bids that would not have their quantity bids fulfilled without the incentive. The additional energy transported is in the form of a blend of natural gas and hydrogen, which adds to the total natural gas transported by the system. Therefore, increasing \(c_m^{CO_2}\) can counter-intuitively lead to more carbon emissions under the specific instance where the system has unmet demands because of the prevailing price structure.

\begin{table}[!h]
\centering
\begin{tabular}{|p{7.2cm}|c|c|c|c|}
 \hline
 Scenario & Baseline & Scenario 1 & Scenario 2 & Scenario 3 \\
 \hline
 \hline
 Total NG Delivered [kg/s] & 727 & 688 & 674 & 679 \\
 Total H\(_2\) Delivered [kg/s] & 45 & 76 & 75 & 75 \\
 Total Energy Delivered [MJ/s] & 38491 & 41223 & 40398 & 40730 \\
 Average Carbon Intensity [kgCO$_2$/MJ] & 0.052 & 0.046 & 0.046 & 0.046 \\
 Total CO\(_2\) [kgCO$_2$/s] & 1998 & 1891 & 1853 & 1868 \\
 Market Revenue \(J_{MR}\) [\$/s] & 550 & 585 & 544 & 544  \\
 Carbon Emissions Incentives \(J_{CEM}\) [\$/s] & 22 & 37 & 36 & 103  \\
 Gas Compression Cost \(J_{GC}\) [\$/s] & 0 & 0.03 & 0.01 & 0.01  \\
 Total Economic Value \(J_{EV}\) [\$/s] & 572 & 622 & 580 & 647  \\
 Total Pass-Through Credits \(D_{PTC}\) [\$/s] & 22 & 37 & 36 & 103  \\
 \hline
\end{tabular}
\caption{Comparison of aggregate optimal outcomes in scenarios for the 40-node network case study. \label{tab:fortynode_results}}
\end{table}

Finally, we can examine the overall economic effect of incentives for decarbonization by hydrogen blending when prevailing energy prices are low in this case study.  In particular, we compare the outcomes in Scenarios 2 and 3 by inspecting the objective function components in Table \ref{tab:fortynode_results}.  While the total amount of hydrogen delivered does not change, the total energy delivered increases from 40398 to 40730 MJ/s, and that increase is in the form of an extra 5 kg/s of natural gas. At the same time, the total economic value produced by the pipeline increases from \$580/s to \$647/s, an increase of \$67/s.  However, note that the total pass-through credits increase from \$37/s to \$103/s, an increase of \$66/s.  It follows that the entire increase in incentives to decarbonize is passed through to customers who as a result only consume more natural gas and \emph{increase} carbon emissions.  Based on this analysis, we argue that it is paramount to conduct engineering economic analyses using physics-based optimization to ensure that any hydrogen blending initiatives and incentives do in fact result in decreases carbon emissions.

\section{Conclusion}  \label{sec:conclusion}

We demonstrated an optimization-based market mechanism for pricing energy transported through a pipeline network that carries mixtures of natural gas and hydrogen and includes decarbonization incentives, and analytically derived the premium paid by consumers that receive hydrogen to lower their carbon intensity.  Our case studies show that when the pipeline is not constrained by hydrogen injection limits, introducing a carbon emissions reduction incentive leads to a decrease in overall natural gas consumption and consequently to a reduction in total carbon dioxide emitted. However, accounting for the effect of  customers whose quantity bids are not fulfilled because of their low price bids, we encounter counter-intuitive behavior of solutions to problem \eqref{eq:optimization_formulation} in which increasing the decarbonization incentive can lead to greater carbon emissions.  The results for the 40-node network extend our previous results \cite{sodwatana2023h2blend,zlotnik2023review} by providing an analysis of optimality conditions and demonstrating the scalability of the computational implementation to a large, complex test case with multiple loops and numerous injection and withdrawal sites. The computational examples also show that optimization-based analysis is critical for making sound decisions about economic policy for blending green hydrogen into existing natural gas pipeline systems.  Moreover, we demonstrate that the problem \eqref{eq:optimization_formulation} can be used by a pipeline market administrator not only to schedule physical and financial pipeline operations, but also as a financial mechanism to collect and distribute decarbonization incentives in a manner that is consistent with customer participation, pipeline network structure, and the physics of energy flow.   Future studies could focus on applying optimization-based market mechanisms similar to that presented here to regional pipeline network models to guide the development of regional hydrogen hubs.


\appendix

\section{Nomenclature}  \label{app:nomenclature}

\begin{table}[!h]
\centering \normalsize
\begin{tabular}{|p{3cm}|p{13cm}|}
  \hline
  \multicolumn{2}{|c|}{Decision Variables} \\ 
  \hline
  \hline
  $P_i$   & Pressure at the ingoing node [Pa] \\
  $P_j$   & Pressure at the outgoing node [Pa] \\
  $V_{ij}$   & Linear combination of squared wave speeds [(m/s$^2$)] \\
  $\phi_{ij}$   & Total mass flow along the pipe [kg/s] \\
  $\gamma_{j}$   & Mass fraction of H$_2$ at node $j$ [-] \\
  $\gamma_{ij}$  & Mass fraction of H$_2$ along the pipe [-] \\
  $s_m^{NG}$   & Mass flow rate of natural gas at injection gNode $m$ [kg/s] \\
  $s_m^{H_2}$   & Mass flow rate of H$_2$ at injection gNode $m$ [kg/s] \\
  $d_m$   & Mass flow rate of delivered blended gas [kg/s] \\
  $\alpha_{ij}$   & Boost ratio of the compressor between suction $i$ and discharge $j$ [-] \\
  \hline
\end{tabular}
\caption{Decision variables. \label{tab:decision_variables}}
\end{table}

\begin{table}[!h]
\centering \normalsize
\begin{tabular}{|p{3cm}|p{13cm}|}
  \hline
  \multicolumn{2}{|c|}{Derived Quantities} \\ 
  \hline
  \hline
  $W_c$   & Amount of work required to compress gas [kw] \\
  $R(\gamma_{j(m)})$   & Calorific value of the blended gas based on H$_2$ mass fraction [MJ/kg] \\
  $E_m$ & Avoided CO$_2$ emissions at gNode $m$ [kgCO$_2$/s] \\
  \hline
\end{tabular}
\caption{Decision variables. \label{tab:derived_quantities}}
\end{table}

\begin{table}[!h]
\centering \normalsize
\begin{tabular}{|p{3cm}|p{13cm}|}
  \hline
  \multicolumn{2}{|c|}{Physical Parameters} \\ 
  \hline
  \hline
  $a_{H_2}$   & Wave speed of H$_2$ [m/s] \\
  $a_{NG}$   & Wave speed of natural gas [m/s] \\
  $a_{0}$   & Geometric mean of wave speeds in gases [m/s] \\
  $\kappa_{H_2}$   & Specific heat ratio of H$_2$ [-] \\
  $\kappa_{NG}$   & Specific heat ratio of natural gas [-] \\
  $\kappa_{nom}$   & Specific heat ratio of the blended gas when $\gamma_{ij} = 0.05$ [-] \\
  $G_{H_2}$   & Specific gravity of H$_2$ [-] \\
  $G_{NG}$   & Specific gravity of natural gas [-] \\
  $G_{nom}$   & Specific gravity of the blended gas when $\gamma_{ij} = 0.05$ [-]  \\
  $R_{H_2}$   & Calorific value of H$_2$ [MJ/kg] \\
  $R_{NG}$   & Calorific value of natural gas [MJ/kg] \\
  $\zeta$    & Approximate ratio of the molecular weight of CO$_2$ to methane [-] \\
  $R$   & Universal gas constant [J/mol/K] \\
  $M_{H_2}$   & Molecular mass of H$_2$ [kg/mol] \\
  $M_{NG}$   & Molecular mass of natural gas [kg/mol] \\
  \hline
\end{tabular}
\caption{Physical parameters. \label{tab:physical_parameters}}
\end{table}

\clearpage

\begin{table}[!t]
\centering \normalsize
\begin{tabular}{|p{3cm}|p{13cm}|}
  \hline
  \multicolumn{2}{|c|}{Operational Parameters} \\ 
  \hline
  \hline
  $f_{ij}$   & Friction factor of pipe [-] \\
  $L_{ij}$   & Length of pipe [m] \\
  $D_{ij}$   & Diameter of pipe [m] \\
  $A_{ij}$   & Area of pipe [m$^2$] \\
  $\sigma_{j}$   & Pressure at the slack node [Pa] \\
  $T$   & Temperature of gas entering the compressor [K] \\
  $\eta$   & Conversion factor for economic cost of compressor work  [\$/kw-s] \\
  $P_j^{min}$   & Minimum operational pressure limit [Pa] \\
  $P_j^{max}$   & Maximum operational pressure limit [Pa] \\
  $\gamma_j^{min}$ & Minimum operational H$_2$ mass fraction [-] \\
  $\gamma_j^{max}$ & Maximum operational H$_2$ mass fraction [-] \\
  $\alpha_{ij}^{max}$   & Maximum compressor boost ratio [-] \\
  $s_m^{max,NG}$   & Natural gas supplier offer quantity at gNode $m$ [kg/s] \\
  $s_m^{max,H_2}$   & Hydrogen supplier offer quantity at gNode $m$ [kg/s] \\
  $g_m^{max}$   & Bid quantity for energy of a flexible customer at gNode $m$ [MJ/s] \\
  $\bar{g}_m$   & Bid quantity for energy of a fixed customer at gNode $m$ [MJ/s] \\
  $c_m^{NG}$   & Natural gas supplier offer price [\$/kg] \\
  $c_m^{H_2}$   & Hydrogen supplier offer price [\$/kg] \\
  $c_m^{d}$   & Customer bid price for energy [\$/MJ] \\
  $c_m^{CO_2}$  & Incentive for avoided CO$_2$ emissions [\$/kgCO$_2$] \\
  \hline
\end{tabular}
\caption{Operational parameters. \label{tab:operational_parameters}}
\end{table}


\section*{Acknowledgements} The authors are grateful to Vitaliy Gyrya, Luke S. Baker, and Aleksandr M. Rudkevich for insightful discussions.  This study was supported by the U.S. Department of Energy's Advanced Grid Modeling (AGM) project ``Dynamical Modeling, Estimation, and Optimal Control of Electrical Grid-Natural Gas Transmission Systems'', as well as LANL Laboratory Directed R\&D project ``Efficient Multi-scale Modeling of Clean Hydrogen Blending in Large Natural Gas Pipelines to Reduce Carbon Emissions''. Research conducted at Los Alamos National Laboratory is done under the auspices of the National Nuclear Security Administration of the U.S. Department of Energy under Contract No. 89233218CNA000001.  Report No. LA-UR-23-32525.





\bibliographystyle{unsrt}  

\bibliography{references}

\end{document}